\documentstyle[amscd,amssymb,verbatim,12pt]{amsart}
\newcommand{\C}{{\Bbb C}}

\newcommand{\ch}{\mbox{\rm ch}}
\newcommand{\codim}{\mbox{\rm codim}}

\newcommand{\cone}{\mbox{\rm cone}}

\newcommand{\Diff}{\mbox{\rm Diff}}

\newcommand{\End}{\mbox{\rm End}}

\newcommand{\Fin}{\mbox{\rm Fin}}
\newcommand{\Func}{\mbox{\rm Func}}

\newcommand{\GL}{\mbox{\rm GL}}
\newcommand{\hA}{\widehat{A}} 
\newcommand{\hcA}{\widehat{\cal A}}

\newcommand{\HH}{\mbox{\rm H}}

\newcommand{\Id}{\mbox{\rm Id}}

\newcommand{\Image}{\mbox{\rm Im}}

\newcommand{\Ind}{\mbox{\rm Ind}}

\newcommand{\Int}{\mbox{\rm int}}

\newcommand{\Ker}{\mbox{\rm Ker}}
\newcommand{\KO}{\mbox{\rm KO}}

\newcommand{\Man}{\mbox{\rm Man}}

\newcommand{\Or}{\mbox{\rm Or}}

\newcommand{\pt}{\mbox{\rm pt}}
\newcommand{\Q}{{\Bbb Q}}
\newcommand{\R}{{\Bbb R}}

\newcommand{\sign}{\mbox{\rm sign}}

\newcommand{\spin}{\mbox{\rm spin}}

\newcommand{\TR}{\mbox{\rm TR}}

\newcommand{\Z}{{\Bbb Z}}
\theoremstyle{plain}
\newtheorem{definition}{Definition}
\newtheorem{assumption}{Assumption}
\newtheorem{lemma}{Lemma}
\newtheorem{theorem}{Theorem}
\newtheorem{proposition}{Proposition}
\newtheorem{corollary}{Corollary}

\newtheorem{conjecture}{Conjecture}

\numberwithin{equation}{section}
\renewcommand{\rm}{\normalshape}
\begin{document}
\title{Signatures and Higher Signatures of $S^1$-Quotients}
\author{John Lott}
\address{Department of Mathematics\\
University of Michigan\\
Ann Arbor, MI  48109-1109\\
USA}
\email{lott@@math.lsa.umich.edu}
\thanks{Research supported by NSF grant DMS-9704633}
\date{July 13, 1998}
\maketitle
\begin{abstract}
We define and study
the signature, $\widehat{A}$-genus and higher signatures of the 
quotient space of
an $S^1$-action on a closed oriented manifold. We give applications
to questions of positive scalar curvature and to an Equivariant
Novikov Conjecture.
\end{abstract}
\section{Introduction} \label{Introduction}

The signature $\sigma(M)$ is a classical 
invariant of a closed oriented manifold $M$. If $M$ is smooth, the
Hirzebruch signature theorem expresses $\sigma(M)$ in terms of the
$L$-class of $M$ \cite{Hirzebruch (1956)}.
It is of interest to extend Hirzebruch's formula to various types of 
singular spaces.
For example, Thom defined the $L$-class of a $PL$-manifold
in a way so that the
Hirzebruch formula still holds \cite{Thom (1958)}.
A large generalization of Thom's result was given by 
Cheeger and Goresky-MacPherson  
\cite{Cheeger (1983),Goresky-MacPherson (1980)}, who defined the signatures
and homology $L$-classes of so-called Witt spaces.

In this paper we define and study the signatures of 
certain singular spaces which
arise in transformation group theory, namely quotients of closed oriented
smooth manifolds $M$ by $S^1$-actions. This class of spaces includes oriented
manifolds-with-boundary, but also contains spaces with much more drastic
singularities. If the group action is semifree, meaning that each
isotropy subgroup is $\{e\}$ or $S^1$, then any point in the quotient space
$S^1 \backslash M$ which is in the singular stratum
has a neighborhood which is homeomorphic to 
$D^k \times \cone(\C P^N)$, for some $k$ and $N$. 
If $N$ is even then the quotient space is not a Witt space. 

Our motivation to study such spaces 
comes from the Equivariant Novikov Conjecture. The usual
Novikov Conjecture hypothesizes that the higher signatures of a closed
oriented manifold are oriented-homotopy invariants.  When one studies
compact group actions, one wants to know what the possible 
equivariant homotopy invariants are. In particular, in view of the importance
of the Novikov conjecture in surgery theory, one wants to know if
there are equivariant higher signatures and an Equivariant Novikov
Conjecture.
There are two candidate Equivariant
Novikov Conjectures, one based on classifying spaces for compact group actions
\cite{Rosenberg-Weinberger (1990)} and one based on classifying spaces for
proper group actions \cite[Section 8]{Baum-Connes-Higson (1994)}. 
We describe these in
detail in Subsection \ref{Equivariant Novikov Conjectures}. In the special 
case of free $S^1$-actions on simply-connected manifolds, 
the first conjecture is false (as was 
pointed out in \cite{Rosenberg-Weinberger (1990)}) and the second conjecture
is true but vacuous.  Since the usual signature of the quotient space of
a free $S^1$-action is an oriented $S^1$-homotopy invariant, there is clearly 
something missing in these conjectures. Hence it is a serious
conceptual problem
to even give a good notion of equivariant higher signatures.

We start out by considering the case when there are no fundamental group
complications.
In Section \ref{Signatures of $S^1$-Quotients}
we define the equivariant signature $\sigma_{S^1}(M) \in \Z$ of an
$S^1$-action. In the special case when $S^1 \backslash M$ is a manifold
(possibly with boundary), $\sigma_{S^1}(M)$ equals the usual signature of
$S^1 \backslash M$.

Note that the fixed-point-set $M^{S^1}$ embeds in
$S^1 \backslash M$.   
Let $\int_{S^1 \backslash M}$ denote integration over
$\left( S^1 \backslash M \right) - M^{S^1}$.
\begin{theorem} \label{th1}
$\sigma_{S^1}(M)$ is an oriented $S^1$-homotopy invariant.
Suppose that the $S^1$-action is semifree.
If $M$ is equipped with an $S^1$-invariant Riemannian metric, give
$\left( S^1 \backslash M \right) - M^{S^1}$ the quotient metric. Then
\begin{equation} \label{eqn1}
\sigma_{S^1}(M) = \int_{S^1 \backslash M} L \left(
T(S^1 \backslash M) \right) +  \eta \left( M^{S^1} \right),
\end{equation}
where $\eta \left( M^{S^1} \right)$ is the Atiyah-Patodi-Singer eta-invariant
of the tangential signature operator on $M^{S^1}$ 
\cite{Atiyah-Patodi-Singer II (1975)}.
\end{theorem}
We also give the extension of 
(\ref{eqn1}) to general $S^1$-actions. If $S^1 \backslash M$ is a Witt space,
we show that $\sigma_{S^1}(M)$ equals the intersection-homology signature
of $S^1 \backslash M$.

In Subsection \ref{Ahat-Genus} we define the
$\widehat{A}$-genus $\widehat{A}_{S^1}(M) \in \Z$ of the quotient space of an 
even semifree
$S^1$-action on a spin manifold $M$. This has applications to questions of
positive scalar curvature. Let us recall the result of
B\'erard-Bergery that if $S^1$ acts freely on a compact manifold $M$ then
$M$ has an $S^1$-invariant metric of positive scalar curvature if and only
if $S^1 \backslash M$ has a metric of positive scalar curvature
\cite[Theorem C]{Berard-Bergery (1983)}. As a consequence,
if $S^1$ acts freely and evenly on a spin manifold $M$ then the usual
$\widehat{A}$-genus of $S^1 \backslash M$ is an obstruction to having an
$S^1$-invariant metric on $M$ of positive scalar curvature. We extend this
to a statement about semifree $S^1$-actions.

\begin{theorem} \label{th1.5}
Suppose that $S^1$ acts semifreely and evenly on a spin manifold $M$. 
If $M$ admits an
$S^1$-invariant metric of positive scalar curvature and $M^{S^1}$ has no
connected components of codimension $2$ in $M$ then $\widehat{A}_{S^1}(M) = 0$.
\end{theorem}

The codimension assumption in Theorem \ref{th1.5} is probably not necessary;
see the remark after the proof of Theorem \ref{th1.5}.

In the case of an odd semifree $S^1$-action on a spin manifold $M$, we
define the $\widehat{A}$-genus $\widehat{\cal A}_{S^1}(M) \in \Z$ 
corresponding to a 
$\spin^c$-structure on the quotient space. We show that $\widehat{A}_{S^1}(M)$
and
$\widehat{\cal A}_{S^1}(M)$ are metric-independent provided that there is
no spectral flow for the Dirac operator on $M^{S^1}$.

In Section \ref{Equivariant Higher Indices} we construct 
equivariant higher signatures for $S^1$-actions, 
using \cite{Lott (1992),Lott (1992a),Lott (1996)}. Let $\Gamma^\prime$
be a finitely-generated discrete group and let 
$\rho : \pi_1(M, m_0) \rightarrow \Gamma^\prime$ be a surjective homomorphism.
Let $o \in \pi_1(M, m_0)$ be the homotopy class of the $S^1$-orbit of 
a basepoint $m_0$
and let $\widehat{\Gamma}$ be the quotient of $\Gamma^\prime$ by the central
cyclic
subgroup generated by $\rho(o)$. The equivariant higher signatures will
involve the group cohomology of $\widehat{\Gamma}$.

In order to construct the equivariant
higher signatures we make a certain ($S^1$-homotopy invariant)
assumption about $M^{S^1}$. Namely, if $F$ is 
a connected component of $M^{S^1}$,
let $\Gamma_F$ be the image of 
$\pi_1(F)$ in $\widehat{\Gamma}$. Let $\overline{\cal D}$ be the
canonical flat $C^*_r \Gamma_F$-bundle on $F$. We assume that
$\HH^* ( F; \overline{\cal D} )$ vanishes in the middle degree if
$F$ is even-dimensional, or in the middle two degrees if $F$ is 
odd-dimensional. We also assume that $\Gamma_F$ is virtually nilpotent or
Gromov-hyperbolic. 

There is a space $\widehat{M}$ on which $\widehat{\Gamma}$ acts properly
and cocompactly, with $\widehat{\Gamma} \backslash \widehat{M} =
S^1 \backslash M$. We needd two pieces of additional data : an $S^1$-invariant
Riemannian metric $g$ on $M$ and a compactly-supported function $H$ on 
$\widehat{M}$ satisfying 
$\sum_{\gamma \in \widehat{\Gamma}} \gamma \cdot H = 1$.
Given $[\tau] \in \HH^k \left( \widehat{\Gamma}; \R \right)$,
represent it by a cocycle $\tau \in Z^k \left(\widehat{\Gamma}; \R \right)$.
There is a corresponding cyclic cocycle $Z_\tau \in ZC^k(\R \widehat{\Gamma})$.
Using $g$, $H$ and $Z_\tau$, we define a closed orbifold form 
$\omega_\tau \in \Omega^k \left(
\left( S^1 \backslash M \right) - M^{S^1} \right)$.
In Subsection 
\ref{``Moral'' Fundamental Group of Quotient}, we use $\omega_\tau$ to
give a differential form proof of
a result of Browder and Hsiang \cite[Theorem 1.1]{Browder-Hsiang (1982)}, 
in the case of $S^1$-actions.
Now suppose that the $S^1$-action is semifree.
Then using the higher eta-invariant
$\widetilde{\eta}$ of \cite{Lott (1992a)}, 
the equivariant higher signature is defined to be
\begin{equation} \label{eqn102.}
\langle \sigma_{S^1}(M), [\tau] \rangle =
\int_{S^1 \backslash M} L(T(S^1 \backslash M)) \wedge \omega_\tau
\: + \: c(k) \:  \langle \widetilde{\eta}(M^{S^1}), Z_{\tau} \rangle \in \R. 
\end{equation}
Here $c(k)$ is a certain nonzero constant.
\begin{theorem} \label{th2}
$\langle \sigma_{S^1}(M), [\tau] \rangle$ is independent of $g$ and $H$. 
\end{theorem}
Thus $\langle \sigma_{S^1}(M), [\tau] \rangle$ is a (smooth) topological
invariant of the $S^1$-action. If $S^1$ acts freely on $M$ then
we recover the Novikov higher signatures of $S^1 \backslash M$ in full
generality.
(Note that the assumptions just involve $M^{S^1}$).
We also give the extensions of (\ref{eqn102.}) 
and
Theorem \ref{th2} to general $S^1$-actions. We conjecture that 
$\langle \sigma_{S^1}(M), [\tau] \rangle$ is an oriented $S^1$-homotopy 
invariant of $M$. In Appendix A we outline a proof
of this when $S^1 \backslash M$ is a manifold-with-boundary whose fundamental
group is virtually nilpotent or Gromov-hyperbolic.
 
I thank Mark Goresky, Matthias Kreck,
Eric Leichtnam, Paolo Piazza, Stephan Stolz and Shmuel Weinberger
for helpful discussions.
 
\section{Signatures of $S^1$-Quotients}
\label{Signatures of $S^1$-Quotients}
\subsection{$S^1$-Homotopy Invariance}
\label{$S^1$-Homotopy Invariance}
Let $G$ be a compact 
Lie group and let $G-\Man$ be the category whose objects are
closed oriented 
smooth manifolds on which $G$ acts on the left
by orientation-preserving diffeomorphisms,
and whose morphisms are smooth orientation-preserving $G$-maps.
If $H$ is a closed subgroup of $G$,
let $M^{H}$ denote the points of $M$ which are fixed by $H$.

The most basic $G$-homotopy invariant information of a $G$-manifold $M$
is the collection of finite sets
$\{\pi_0(M^H)\}$. To organize this information coherently, 
let $\Or_G$ be the orbit category of $G$,
with objects given by $G$-homogeneous spaces $G/H$, $H$ closed, and 
morphisms given by $G$-maps.  Let $\Fin$ be the category whose objects are
isomorphism classes of finite sets and whose morphisms 
are set maps.
Then there is a functor $F : G-\Man \rightarrow \Func(\Or_G^{op}, \Fin)$ 
where $F(M) \in \Func(\Or_G^{op}, \Fin)$ 
sends  $G/H \in \Or_G^{op}$ to $\{\pi_0(M^H)\}$.
Given $\mu \in \Func(\Or_G^{op}, \Fin)$, the set of  $G$-manifolds 
$M$ such that $F(M) = \mu$ is
closed under $G$-homotopy equivalence. For example, the notion of an action
being free or semifree is $G$-homotopy invariant.

We now restrict to the case $G = S^1$.
Suppose that $M$ has dimension $4k+1$.
Let $X$ be the vector field on $M$ which generates the $S^1$-action. Let
$i_X : \Omega^*(M) \rightarrow \Omega^{*-1}(M)$ be interior multiplication
by $X$ and let ${\cal L}_X : \Omega^*(M) \rightarrow \Omega^{*}(M)$ be Lie
differentiation by $X$. Let $e : M^{S^1} \rightarrow M$ be the
inclusion of the fixed-point-set. 

\begin{definition} \label{def1}
Define the basic forms on $M$ and $\left(M, M^{S^1} \right)$ by
\begin{align} \label{eqn2}
\Omega^{*,basic}(M) & = \{ \omega \in \Omega^*(M) : i_X \omega = 
{\cal L}_X \omega = 0 \}, \\
\Omega^{*,basic} \left( M, M^{S^1} \right) & = 
\{ \omega \in \Omega^{*, basic}(M) : e^* \omega = 0 \}. \notag
\end{align}
Let $\Omega^{*,basic}_c \left (M - M^{S^1} \right)$ be
the complex of compactly-supported basic forms on $M - M^{S^1}$.
Let $\HH^{*, basic}(M)$, $\HH^{*, basic} \left(M, M^{S^1} \right)$ and
$\HH^{*,basic}_c \left (M - M^{S^1} \right)$ 
be the corresponding cohomology groups.
\end{definition}

\begin{proposition} \label{prop1}
$\HH^{*, basic}(M) \cong \HH^*(S^1 \backslash M; \R)$ and
\begin{equation} \label{eqn3}
\HH^{*, basic} \left(M, M^{S^1} \right) \cong
\HH^{*,basic}_c \left (M - M^{S^1} \right) \cong \HH^* \left( 
S^1 \backslash M, M^{S^1}; \R \right).
\end{equation}
\end{proposition}
\begin{pf}
The fact that $\HH^{*, basic}(M) \cong \HH^*(S^1 \backslash M; \R)$ was
proven in \cite{Koszul (1953)}. Let us briefly recall the proof.
By a Mayer-Vietoris argument, we can reduce to the case when  $H$ is
a subgroup of $S^1$, $V$ is a representation space of $H$, $DV$ is the unit
ball in $V$ and
$M = S^1 \times_H DV$. By a product formula, we 
can also reduce to the case when $V$ has no trivial 
subrepresentations. Let $SV$ be the unit sphere in $V$. Then 
$S^1 \backslash M = H \backslash DV$ is a cone over 
the orbifold $H \backslash SV$.
A Poincar\'e lemma gives 
\begin{equation} \label{eqn4}
\HH^{*, basic}(M) \cong \HH^*(\pt; \R) \cong \HH^*(S^1 \backslash M; \R),
\end{equation}
which proves the claim.

We now do a similar argument for $\HH^{*, basic} \left(M, M^{S^1} \right)$.
We can reduce to the case $M = S^1 \times_H DV$ as above.
The Poincar\'e lemma gives 
\begin{align} \label{eqn5}
\HH^{*, basic} \left(M, M^{S^1} \right) & = 
\HH^{*, basic} \left(S^1 \times_H DV, \pt. \right) =
0 \\
& = \HH^*(H \backslash 
DV, \pt; \R) = \HH^* \left(S^1 \backslash M, M^{S^1}; \R \right). \notag
\end{align}

Finally, there is an obvious chain inclusion
$\Omega^{*,basic}_c \left (M - M^{S^1} \right) \rightarrow
\Omega^{*,basic} \left( M, M^{S^1} \right)$. Using the Poincar\'e lemma,
one can construct a homotopy inverse
$\Omega^{*,basic} \left( M, M^{S^1} \right) \rightarrow
\Omega^{*,basic}_c \left (M - M^{S^1} \right)$. The proposition follows.
\end{pf}
\noindent
{\bf Example : } Let $X^{4k}$ be a compact manifold-with-boundary.
Let $M$ be the manifold obtained by spinning $X$. That is,
$M = \partial (D^2 \times X) = \left( D^2 \times \partial X \right) 
\cup_{S^1 \times \partial X} \left( S^1 \times X \right)$, with the
induced $S^1$-action. Then
$\HH^{*,basic} \left (M, M^{S^1} \right)
\cong \HH^{*}(X, \partial X; \R) \cong \HH^{*}_c(\Int (X); \R)$. \\

\begin{proposition} \label{prop2}
$\HH^{*,basic} \left (M, M^{S^1} \right)$ is an 
$S^1$-homotopy invariant.
\end{proposition}
\begin{pf}
Let $f : M \rightarrow N$ be an $S^1$-homotopy equivalence, with 
$S^1$-homotopy inverse $g : N \rightarrow M$. Then
$f \left( M^{S^1} \right) \subset N^{S^1}$.
Hence the pullback
$f^* : \Omega^{*,basic} \left (N,  N^{S^1} \right) \rightarrow
\Omega^{*,basic} \left (M, M^{S^1} \right)$ is well-defined.
Let $F : [0,1] \times M \rightarrow M$ be an $S^1$-homotopy from the identity
to $g \circ f$. Then there is a pullback 
$F^* : \Omega^{*,basic} \left (M, M^{S^1} \right) \rightarrow
\Omega^*([0,1]) \: {\otimes} \: 
\Omega^{*,basic} \left (M, M^{S^1} \right)$. One can construct
the cochain-homotopy equivalence between 
$\Omega^{*,basic} \left (M, M^{S^1} \right)$ and 
$\Omega^{*,basic} \left (N, N^{S^1} \right)$ by the standard argument.
\end{pf}
{\bf Remark : } It is not 
surprising that $\HH^{*,basic} \left (M, M^{S^1} \right)$ is unchanged
by an $S^1$-isovariant homotopy equivalence, as this would correspond
to a stratum-preserving homotopy equivalence between $M/S^1$ and $N/S^1$.
However, it is perhaps less obvious that it is unchanged by an
$S^1$-equivariant homotopy equivalence.\\

In the rest of the paper, we will deal with
$\HH^{*,basic}_c \left (M - M^{S^1} \right)$ instead of the equivalent
$\HH^{*, basic} \left(M, M^{S^1} \right)$.
Give $M$ an $S^1$-invariant Riemannian metric.
Let $X^* \in \Omega^1(M)$ be the dual $1$-form to $X$, using the Riemannian
metric.  Define $\eta \in \Omega^1 \left( M - M^{S^1} \right)$ by
$\eta = X^*/|X|^2$.
\begin{proposition} \label{prop3}
$d\eta$ is a basic $2$-form on $M - M^{S^1}$.
\end{proposition}
\begin{pf}
By construction, ${\cal L}_X \eta = 0$ and hence ${\cal L}_X \: d \eta = 0$.
Also by construction, $i_X \eta = 1$. Hence $i_x \: d \eta = 
{\cal L}_X \eta - d \: i_X \eta = 0$.
\end{pf}

\begin{proposition} \label{prop4}
If $\sigma \in \Omega^{4k-1,basic}_c \left (M - M^{S^1} \right)$ then
$\int_M \eta \wedge d\sigma = 0$.
\end{proposition}
\begin{pf}
We have
\begin{equation} \label{eqn6}
\int_M \eta \wedge d\sigma = \int_M d\eta \wedge \sigma.
\end{equation}
As $d\eta$ and $\sigma$ are basic, $d\eta \wedge \sigma$ is basic
and so the $(4k+1)$-form $d\eta \wedge \sigma$ vanishes
in $\Omega^{4k+1}_c \left( M - M^{S^1} \right)$. 
\end{pf}

\begin{definition}
The $S^1$-fundamental class of $M$ is the map
$\tau : \HH^{4k,basic}_c \left (M - M^{S^1} \right) \rightarrow \R$ given by
$\tau(\omega) = \int_M \eta \wedge \omega$. \\
\end{definition}

By Proposition \ref{prop4}, the $S^1$-fundamental class is well-defined. \\
\begin{proposition} \label{prop5}
The $S^1$-fundamental class of $M$ is independent of the choice of
Riemannian metric.
\end{proposition}
\begin{pf}
Let $\eta_1$ and $\eta_2$ be the $1$-forms coming from two 
Riemannian metrics. Then $\eta_1 - \eta_2$ is basic. Hence
$\int_M (\eta_1 - \eta_2) \wedge \omega = 0$ for any
$\omega \in \Omega^{4k,basic}_c \left (M - M^{S^1} \right)$.
\end{pf}

\begin{definition}
The intersection form on $\Omega^{2k,basic}_c \left (M - M^{S^1} \right)$ is
\begin{equation} \label{eqn7}
(\omega_1, \omega_2) = \int_M \eta \wedge \omega_1 \wedge \omega_2.
\end{equation}
\end{definition}

Clearly $(\cdot, \cdot)$ is symmetric. By Proposition \ref{prop4}, it 
extends to a bilinear form on $\HH^{2k,basic}_c \left (M - M^{S^1} \right)$.
 
\begin{definition}
$\sigma_{S^1}(M)$ is the signature of $(\cdot, \cdot)$. That is, if
the symmetric form
$(\cdot, \cdot)$ is diagonalized then $\sigma_{S^1}(M)$ is (the number of 
positive eigenvalues) minus (the number of negative eigenvalues).
\end{definition}
\noindent
{\bf Remark : } 
The symmetric form $(\cdot, \cdot)$ on
$\HH^{2k,basic}_c \left (M - M^{S^1} \right)$ 
may be degenerate. For example,
if $M$ comes from spinning an oriented compact manifold-with-boundary $X$ 
then the intersection form 
$(\cdot, \cdot)$ on $\HH^{2k,basic}_c \left (M - M^{S^1} \right)$ is the 
same as that on $\HH^{2k}(X, \partial X)$, 
which may be degenerate. In this case,
$\sigma_{S^1}(M) = \sigma(X)$.\\
\begin{proposition} \label{prop6}
If $f : M \rightarrow N$ is an orientation-preserving $S^1$-homotopy
equivalence then $\sigma_{S^1}(M) = \sigma_{S^1}(N)$.
\end{proposition}
\begin{pf}
It suffices to show that the $S^1$-fundamental class on $M$ pushes forward
to the $S^1$-fundamental class on $N$.
Let $\eta_N$ be the $1$-form constructed from a Riemannian metric on $N$.
For $\omega \in \Omega^{4k,basic}_c \left (N - N^{S^1} \right)$, we have
\begin{equation} \label{eqn8}
\int_N \eta_N \wedge \omega = \int_M f^* \eta_N \wedge f^* \omega.
\end{equation}
However, $f^* \eta_N - \eta_M$ is a basic $1$-form on $M$ and so
\begin{equation} \label{eqn9}
\int_M  \left( f^* \eta_N - \eta_M \right) \wedge f^* \omega = 0.
\end{equation} 
The proposition follows.
\end{pf}

\subsection{Fixed-point-free Actions}
\label{Fixed-point-free Actions}
Let $S^1$ act effectively
on $M$ without fixed points. Then $S^1\backslash M$ is an
oriented orbifold.  If $M$ has an $S^1$-invariant Riemannian metric then
$S^1\backslash M$ is a Riemannian orbifold. To write the
formula for $\sigma_{S^1}(M)$, we first describe a certain set of 
suborbifolds ${\cal O}$
of $S^1 \backslash M$. We construct these suborbifolds by describing their
intersections with orbifold coordinate charts in $S^1 \backslash M$; the
suborbifolds can then be defined by patching together these intersections. 
Given $x \in S^1 \backslash M$, let $\Gamma$ 
be a finite group and let $U \subset \R^n$ be a 
domain with a $\Gamma$-action such that $(\Gamma, U)$ 
is an orbifold coordinate chart for $S^1 \backslash M$ around $x$. 
In particular, $\Gamma \backslash U$ can be identified with a 
neighborhood of $x$. Put 
\begin{equation} \label{eqn10}
\widehat{U} = \{ (g,u) \in \Gamma \times U : g u = g\}.
\end{equation} 
Define a $\Gamma$-action on
$\widehat{U}$ by $\gamma \cdot (g, u) = (\gamma g \gamma^{-1}, g u)$.
Let $\pi : \widehat{U} \rightarrow \Gamma \backslash \widehat{U}$ be the
quotient map.
Let $\langle \Gamma \rangle$ denote the set of conjugacy classes of
$\Gamma$. There are projection maps $p_1 : \Gamma \backslash \widehat{U}
\rightarrow \langle \Gamma \rangle$ and $p_2 : \Gamma \backslash \widehat{U}
\rightarrow \Gamma \backslash U$. Then the intersections of the
${\cal O}$'s with $\Gamma \backslash U$ are
$\{p_2 p_1^{-1}(\langle g \rangle) \}_{\langle g \rangle \in
\langle \Gamma \rangle}$. 

Define the multiplicity $m_x \in \Z^+$ of $x$ by
$\frac{1}{m_x} = 
\frac{|(p_2 \circ \pi)^{-1}(x)|}{|\Gamma|}$. This is independent
of the choice of orbifold coordinate chart. 

The Atiyah-Singer equivariant $L$-class $L(g) \in \Omega^{even}(U^g,
o(TU^g))$ \cite{Atiyah-Singer III (1968)} 
is the pullback of a differential form 
$L(\langle g \rangle)$
on the image of $U^g$ in $\Gamma \backslash U$. 
Given a suborbifold ${\cal O}$, define 
$L({\cal O}) \in \Omega^{even}({\cal O}, o(T{\cal O}))$ by
\begin{equation} \label{eqn11}
L({\cal O}) \big|_{{\cal O} \cap (\Gamma \backslash U)} = 
\sum_{\langle g \rangle : p_2 p_1^{-1}(\langle g \rangle) =
{\cal O} \cap (\Gamma \backslash U)} L( \langle g \rangle).
\end{equation}

If ${\cal O}$ is one of the
suborbifolds then $m_x$ is constant on the regular
part of ${\cal O}$ and so we can define the multiplicity 
$m_{\cal O} \in \Z^+$ of ${\cal O}$. From \cite{Kawasaki (1978)}, 
it follows that
\begin{equation} \label{eqn12}
\sigma_{S^1}(M) = \sum_{\cal O} \frac{1}{m_{\cal O}}
\int_{\cal O} L({\cal O}).
\end{equation}

By definition, $\sigma_{S^1}(M) \in \Z$. In fact, it equals
the signature of $S^1 \backslash M$ as a rational homology manifold.
In the orbifold world it may be more natural to consider the
${\Q}$-valued orbifold signature 
$\int_{S^1 \backslash M} L(T(S^1 \backslash M))$. 
However, this is definitely a different object and is a single term in
(\ref{eqn12}). \\ \\
{\bf Remark :} In the case of fixed-point-free actions, 
$\sigma_{S^1}(M)$ comes from the index of a signature operator which is
transversally elliptic in the sense of \cite{Atiyah (1974)}. This
transversally elliptic signature operator only exists in the
fixed-point-free case.

\subsection{Semifree Actions}
\label{Semifree Actions}
Suppose that $S^1$ acts effectively and 
semifreely on $M$. Let $(S^1 \backslash M) - M^{S^1}$ have the
quotient Riemannian metric.  We write $L \left(
T(S^1 \backslash M) \right)$ for the $L$-form and
$\int_{S^1 \backslash M} L \left(
T(S^1 \backslash M) \right)$ for its integral over 
$(S^1 \backslash M) - M^{S^1}$.
\begin{theorem} \label{prop7}
\begin{equation} \label{eqn13}
\sigma_{S^1}(M) = \int_{S^1 \backslash M} L \left(
T(S^1 \backslash M) \right) +  \eta \left( M^{S^1} \right).
\end{equation}
\end{theorem}
\begin{pf}
Let $F$ be a connected component of the
fixed-point-set $M^{S^1}$. It is an oriented odd-dimensional manifold, 
say of dimension $4k-2N+1$. Let $NF$ be the normal bundle of $F$ in $M$. 
It has an $S^1$-action by orthogonal automorphisms, which is
free on $NF - F$. Furthermore, the disk bundle
$DNF$ is $S^1$-diffeomorphic to a neighborhood of $F$ in $M$. 
Let $SNF$ be the sphere bundle of $NF$. Then $S^1 \backslash
SNF$ is the total space of a Riemannian
fiber bundle ${\cal F}$ over $F$ whose fibers $Z$ are 
copies of $\C P^N$. The quotient space $S^1 \backslash
DNF$ is homeomorphic to the mapping cylinder of the projection $\pi :
{\cal F} \rightarrow F$.  

Let us first pretend that for each $F$, a neighborhood of $F$ in 
$M$ is $S^1$-isometric to $DNF$. For simplicity,
we suppose that there is only one connected component $F$ of $M^{S^1}$;
the general case is similar.
For $r > 0$, let $N_r(F)$ 
be the $r$-neighborhood of $F$ in $S^1 \backslash M$.
Then for small $r$, $\sigma_{S^1}(M) = 
\sigma \left( (S^1 \backslash M) - N_r(F) \right)$. By the Atiyah-Patodi-Singer
theorem,
\begin{equation} \label{eqn14}
\sigma \left( (S^1 \backslash M) - N_r(F) \right) = 
\int_{(S^1 \backslash M) - N_r(F)} L \left( T(S^1 \backslash M) \right)
+ \int_{\partial N_r(F)} \widetilde{L}(\partial N_r(F)) + \eta(
\partial N_r(F)),
\end{equation} 
where $\widetilde{L}(\partial N_r(F))$ 
is a local expression on $\partial N_r(F)$ which
involves the second fundamental form and the curvature tensor
\cite[Section 3.10]{Gilkey (1995)} and we give $\partial N_r(F)$ the
orientation induced from that of $N_r(F)$. We will compute
the limit of the right-hand-side of (\ref{eqn14}) as $r \rightarrow 0$.

We use the notation of \cite[Section III(c)]{Bismut-Lott (1995)} 
to describe the geometry
of the fiber bundle ${\cal F}$. In particular, the second fundamental form
of the fibers and the curvature of the fiber bundle are parts of
the connection $1$-form component 
$\omega^i_{\: \alpha} = \omega^i_{\: \alpha j} \: \tau^j \: + \:
\omega^i_{\: \alpha \beta} \: \tau^\beta$. Let $R^{TZ}$ be the curvature of
the Bismut connection on $TZ$. Define $\widehat{\Omega} \in
\Omega^2([0,1] \times {\cal F}) \otimes
\End(TZ \oplus \R)$, a skew-symmetric matrix of $2$-forms, by
\begin{align} \label{eqn15}
\widehat{\Omega}^i_{\: j} & = (R^{TZ})^i_{\: j} - t^2 \tau^i \wedge \tau^j, \\
\widehat{\Omega}^i_{\: r} & = dt \wedge \tau^i - t \: 
\omega^i_{\: \alpha} \wedge \tau^\alpha. \notag
\end{align}

\begin{definition}
The transgressed $L$-class, $\widehat{L}(F) \in \Omega^{odd}(F)$, is given by
\begin{equation} \label{eqn16}
\widehat{L}(F) = \int_Z \int_0^1 L \left( \widehat{\Omega} \right).
\end{equation}
\end{definition}

We first compute the curvature of $S^1 \backslash DNF$ 
in terms of the geometric invariants
of the fiber bundle ${\cal F}$. Let $\{\tau^i\}_{i=1}^{dim(Z)},
\{\tau^\alpha\}_{\alpha = 1}^{dim(F)}$ be a local orthonormal basis of 
$1$-forms on ${\cal F}$ as in 
\cite[Section III(c)]{Bismut-Lott (1995)}. 
Then a local orthonormal basis of $1$-forms on
$S^1 \backslash DNF$ is given by
\begin{align} \label{eqn17}
\widehat{\tau}^r & = dr, \\
\widehat{\tau}^i & = r \: \tau^i, \notag \\
\widehat{\tau}^\alpha & = \tau^\alpha. \notag
\end{align}
Let $\omega^I_{\: J}$ be the connection matrix of ${\cal F}$.
The structure equations $0 = d \widehat{\tau}^I + 
\widehat{\omega}^I_{\: J} \wedge
\widehat{\tau}^J$ give the connection matrix
$\widehat{\omega}$ of $S^1 \backslash DNF$ to be
\begin{align} \label{eqn18}
\widehat{\omega}^i_{\: j} & = {\omega}^i_{\: j}, \\
\widehat{\omega}^i_{\: r} & = \tau^i, \notag \\
\widehat{\omega}^i_{\: \alpha} & = r \: {\omega}^i_{\: \alpha}, \notag \\
\widehat{\omega}^\alpha_{\: \beta}
& = r^2 \: {\omega}^\alpha_{\: \beta i} \: \tau^i 
\: + \: {\omega}^\alpha_{\: \beta \gamma} \: \tau^\gamma, \notag \\
\widehat{\omega}^\alpha_{\: r}
& = 0. \notag
\end{align}
In the limit when $r \rightarrow 0$, the curvature matrix of $S^1 \backslash
DNF$ becomes
\begin{align} \label{eqn19}
\widehat{\Omega}^i_{\: j} & = (R^{TZ})^i_{\: j} 
\: - \: \tau^i \wedge \tau^j, \\
\widehat{\Omega}^i_{\: r} & = - \omega^i_{\: \alpha} \wedge \tau^i, \notag \\
\widehat{\Omega}^i_{\: \alpha} & = dr \wedge 
{\omega}^i_{\: \alpha}, \notag \\
\widehat{\Omega}^\alpha_{\: \beta} & = (R^{TF})^\alpha_{\: \beta}, \notag \\
\widehat{\Omega}^\alpha_{\: r}
& = 0. \notag
\end{align}
(By way of illustration, 
let us check the $\widehat{\Omega}^i_{\: r}$-term explicitly.
We have
\begin{align} \label{eqn20}
\widehat{\Omega}^i_{\: r} & = d \widehat{\omega}^i_{\: r} \: + \:
\widehat{\omega}^i_{\: j} \wedge \widehat{\omega}^j_{\: r} \\
& = d \tau^i \: + \:
{\omega}^i_{\: j} \wedge \tau^j \notag \\
& = - \: \omega^i_{\: j} \wedge \tau^j
\: - \: \omega^i_{\: \alpha} \wedge \tau^\alpha
\: + \:
{\omega}^i_{\: j} \wedge \tau^j \notag \\
& =  - \: \omega^i_{\: \alpha} \wedge \tau^\alpha.) \notag
\end{align}

It is now clear that 
\begin{equation} \label{eqn21}
\int_{S^1 \backslash M} 
L \left( T \left(S^1 \backslash M \right) \right) =
\lim_{r \rightarrow 0} 
\int_{(S^1 \backslash M) - N_r(F)} L \left( T(S^1 \backslash M) \right)
\end{equation}
exists.

Restricted to $\partial N_r(F)$, as $r \rightarrow 0$
the curvature matrix has nonzero entries
\begin{align} \label{eqn22}
\widehat{\Omega}^i_{\: j} & = (R^{TZ})^i_{\: j} 
\: - \: \tau^i \wedge \tau^j, \\
\widehat{\Omega}^i_{\: r} & = - \omega^i_{\: \alpha} \wedge \tau^i, \notag \\
\widehat{\Omega}^\alpha_{\: \beta} & = (R^{TF})^\alpha_{\: \beta}. \notag
\end{align}
The second fundamental form of $\partial N_r(F)$ enters in the connection
matrix element
\begin{equation} \label{eqn23}
\widehat{\omega}^i_{\: r} = \tau^i = \frac{1}{r} \: \widehat{\tau}^i.
\end{equation}
That is, with respect to the orthogonal decomposition $T(\partial N_r(F)) =
TZ \oplus \pi^* TF$, the shape operator of $\partial N_r(F)$ is
$\frac{1}{r} \: 
\Id_{TZ} \oplus 0$. 

To compute $\widetilde{L}(\partial N_r(F))$ for
small $r$, we construct a $1$-parameter family of connections which
interpolate between the Riemannian connection of $S^1 \backslash DNF$ 
(pulled back to 
$\partial N_r(F)$) and the Riemannian connection of a product metric, at
least as $r \rightarrow 0$.
For $t \in [0,1]$, put 
\begin{align} \label{eqn24}
\widehat{\omega}^i_{\: j}(t) & = {\omega}^i_{\: j}, \\
\widehat{\omega}^i_{\: r}(t) & = t \: \tau^i, \notag \\
\widehat{\omega}^\alpha_{\: \beta}(t) & = 
{\omega}^\alpha_{\: \beta \gamma} \: \tau^\gamma. \notag
\end{align}
Then $\widehat{\omega}(0)$ is the product connection $1$-form and
by (\ref{eqn18}),
$\widehat{\omega}(1)$ is the limit of the pullback connection $1$-form on 
$\partial N_r(F)$ as $r \rightarrow 0$.
The curvature of (\ref{eqn24}) on $[0,1] \times {\cal F}$ has 
nonzero components
\begin{align} \label{eqn25}
\widehat{\Omega}^i_{\: j}(t) & = (R^{TZ})^i_{\: j} 
- t^2 \tau^i \wedge \tau^j, \\
\widehat{\Omega}^i_{\: r}(t) & = dt \wedge \tau^i - t \: 
\omega^i_{\: \alpha} \wedge \tau^\alpha, \notag \\
\widehat{\Omega}^\alpha_{\: \beta}(t) & = 
(R^{TF})^\alpha_{\: \beta}. \notag
\end{align}
Then 
\begin{equation} \label{eqn26}
\lim_{r \rightarrow 0} \int_{\partial N_r(F)} 
\widetilde{L}(\partial N_r(F)) =
\int_{[0,1] \times {\cal F}} L(\widehat{\Omega}(t)) =
\int_F \widehat{L}(F) \wedge L(TF).
\end{equation}

Now the fiber bundle $S^1 \backslash SNF$ 
is associated to a principal bundle $P$ over $F$ with compact structure group.
Hence $\widehat{L}(F)$ can be computed by equivariant methods 
\cite[Section 7.6]{Berline-Getzler-Vergne (1992)}.
Such a calculation will necessarily give it as a 
polynomial in the curvature form of $P$, and in particular as an even form
on $F$. 
However, by parity reasons, $\widehat{L}(F)$ is an odd form on $F$. Thus
\begin{equation} \label{eqn27}
\lim_{r \rightarrow 0} \int_{\partial N_r(F)} \widetilde{L}(\partial N_r(F))
= 0.
\end{equation}

From \cite{Dai (1991)},
\begin{equation} \label{eqn28}
\lim_{r \rightarrow 0} \eta(\partial N_r(F)) = \int_F \widetilde{\eta}
\wedge L(TF) + \eta(F; \Ind(D_Z)) + \tau_F,
\end{equation}
where $\tau_F$ is a signature correction term \cite[p. 268]{Dai (1991)}.
Again, we can compute $\widetilde{\eta}$ by equivariant methods to obtain
an even form on $F$, while by parity reasons 
$\widetilde{\eta}$ is an odd form. Thus 
\begin{equation} \label{eqn29}
\int_F \widetilde{\eta} \wedge L(TF) = 0.
\end{equation}
Next, the index bundle $\Ind(D_Z)$ on $F$ 
is the difference of the vector bundles
$\HH^{N}_+(Z)$ and $\HH^{N}_-(Z)$ of self-dual and anti-self-dual
cohomology groups. As $Z = \C P^N$,
$\HH^{N}_\pm(Z; \R)$ vanishes unless $N$ is even, in which case 
$\HH^{N}_+(Z; \R) = \R$ and $\HH^{N}_-(Z; \R) = 0$.
Then $\HH^{N}_+(Z)$ is a trivial real line
bundle on $F$ with a flat Euclidean metric. Thus
$\eta(F; \Ind(D_Z)) = \eta(F)$.

Finally, 
the Leray-Hirsch theorem implies that the Leray-Serre spectral
sequence for $\HH^*({\cal F}; \R)$ degenerates at the $E_2$-term
\cite[p. 170, 270]{Bott-Tu (1982)}. Hence $\tau_F = 0$.

This proves the proposition if a neighborhood of 
$F$ in $M$ is $S^1$-isometric to $DNF$.
If a neighborhood of $F$ in $M$ is not $S^1$-isometric to
$DNF$, nevertheless as one approaches $F$ the Riemannian
metric on $M$ is better and better approximated by
that of $DNF$. The above calculations will still be valid
in this limit.
\end{pf}
{\bf Example :} If $M$ is obtained by spinning a compact oriented 
manifold-with-boundary $X$ then $M^{S^1} = - \partial X$, when one takes
orientations into account.  The boundary $\partial X$ in $X = S^1 \backslash M$
is totally geodesic.
In this case, Theorem \ref{prop7} reduces to
the Atiyah-Patodi-Singer formula for $\sigma(X)$.

\begin{proposition} \label{prop8} Let $W$ be a semifree $S^1$-cobordism between
$M_1$ and $M_2$. Then
\begin{equation} \label{eqn30}
\sigma_{S^1} \left(M^1 \right) - \sigma_{S^1}\left(M^2 \right) = 
\sigma \left( W^{S^1} \right).
\end{equation}
\end{proposition}
\begin{pf}
We take $\partial W = M_1 \cup (-M_2)$.
Let $N_r \left( W^{S^1} \right)$ be the $r$-neighborhood of $W^{S^1}$ in 
$S^1 \backslash W$. Then for $r$ small,
$\left( \left( S^1 \backslash W \right) - N_r \left( W^{S^1} \right),
\partial N_r \left( W^{S^1} \right) \right)$ is a cobordism of manifold pairs
from $\left( \left( S^1 \backslash M_1 \right) - N_r \left( M_1^{S^1} \right),
\partial N_r \left( M_1^{S^1} \right) \right)$ to
$\left( \left( S^1 \backslash M_2 \right) - N_r \left( M_2^{S^1} \right),
\partial N_r \left( M_2^{S^1} \right) \right)$. 
Hence $\partial N_r \left( W^{S^1} \right) \cup
\left(  \left( S^1 \backslash M_1 \right) - N_r \left( M_1^{S^1} \right) 
\right)
\cup -
\left( \left( S^1 \backslash M_2 \right) - N_r \left( M_2^{S^1} \right) 
\right)$ is an oriented 
boundary, where $\partial N_r \left( W^{S^1} \right)$ has
the boundary orientation coming from $N_r \left( W^{S^1} \right)$. 
Giving it the other orientation, we obtain
\begin{equation} \label{eqn31}
\sigma_{S^1} \left(M^1 \right) - \sigma_{S^1}\left(M^2 \right) = 
\sigma \left( \partial N_r \left( W^{S^1} \right) \right).
\end{equation}
Now $\partial N_r \left( W^{S^1} \right)$ is the total space of a fiber
bundle with fiber $Z$, a complex projective space, and base
$W^{S^1}$. By the same calculation as at the end of the proof of
Theorem \ref{prop7}, the boundary fibration over $\partial W^{S^1} =
M_1^{S^1} \cup (- M_2^{S^1})$ has vanishing signature correction $\tau$.
Then by \cite[Theorem 0.4b, p. 315]{Dai (1991)},
\begin{equation} \label{eqn32}
\sigma \left(\partial N_r \left( W^{S^1} \right) \right) = \sigma(Z) \cdot
\sigma \left( W^{S^1} \right) = \sigma \left( W^{S^1} \right).
\end{equation}
The proposition follows.
\end{pf}

By way of comparison, the $S^1$-semifree cobordism-invariant information of $M$
essentially consists of the cobordism classes of the components of $M^{S^1}$, 
listed by dimension,
along with their normal data \cite{Uchida (1970)}.

If the codimension of $M^{S^1}$ in $M$ is divisible by four then
$S^1 \backslash M$ is a Witt space in the sense of \cite{Siegel (1983)}.
Hence it has an $L$-class in $\HH_*(S^1 \backslash M; \Q)$.
Also, as $\dim \left( M^{S^1} \right) \equiv 1 \mod{4}$, the eta-invariant
of $M^{S^1}$ vanishes. We use the differential-form description of the
homology of Witt spaces given in  
\cite[Section 4]{Bismut-Cheeger (1991)}.
\begin{proposition} \label{prop9}
In this case, the homology $L$-class of $S^1 \backslash M$ is represented
by $L(T(S^1 \backslash M))$.
\end{proposition}
\begin{pf}
We can deform the metric of $S^1 \backslash M$ to make it strictly conical
in a neighborhood of $M^{S^1}$.
By \cite[Theorem 5.7]{Bismut-Cheeger (1991)}, 
the homology $L$-class is represented
by the pair of forms $\left( L(T(S^1 \backslash M)), L(TM^{S^1}) \wedge 
\widetilde{\eta} \right)$, where $\widetilde{\eta}$ is the eta-form of
the $\C P^N$-bundle over $M^{S^1}$. By the method of
proof of Theorem \ref{prop7},
$\widetilde{\eta} = 0$.
\end{pf}
\begin{corollary}
In this case, $\sigma_{S^1}(M) = \int_{S^1 \backslash M}
L(T(S^1 \backslash M))$ equals the intersection homology signature
of $S^1 \backslash M$.
\end{corollary}

One can give a more direct proof of the corollary.  For small $r > 0$, let 
$N_r \left(M^{S^1} \right)$ be the $r$-tubular neighborhood of $M^{S^1}$
in $S^1 \backslash M$.
As in \cite[Proposition 3.1]{Siegel (1983)}, there is a Witt cobordism
which pinches $\partial N_r \left(M^{S^1} \right)$ to a point. Letting $X_1$ be
the coning of $(S^1 \backslash M) - N_r \left(M^{S^1} \right)$ and $X_2$
be the coning of $\overline{N_r \left(M^{S^1} \right)}$, it follows that
\begin{equation} \label{eqn33}
\sigma(S^1 \backslash M) = \sigma(X_1) + \sigma(X_2),
\end{equation} 
where $\sigma$ denote the intersection homology signature.  Now
$\sigma(X_1) = \sigma_{S^1}(M)$. Let $X_3$ be the mapping cylinder of
the projection $\overline{N_r \left(M^{S^1} \right)} \rightarrow M^{S^1}$.
A further Witt cobordism shows that $\sigma(X_2) = \sigma(X_3)$. 
It is well-known that the signature of the total space of an oriented fiber
bundle vanishes if the fiber and base have odd dimension.  One can extend
this fact to the fibration $X_3 \rightarrow M^{S^1}$, whose fiber is a
Witt space, as in \cite[p. 545-546]{Cappell-Shaneson (1991)}.
(Strictly speaking, \cite{Cappell-Shaneson (1991)} deals with the more
interesting case of even-dimensional fiber and base.) The corollary follows.

We expect that for a general semifree effective
$S^1$-action, $\sigma_{S^1}(M)$ will be the
signature of the intersection pairing on the image of the (lower middle
perversity) middle-dimensional intersection homology in the (upper middle
perversity) middle-dimensional intersection homology.

\subsection{$\widehat{A}$-Genus}
\label{Ahat-Genus}

We wish to construct an analog of the $\widehat{A}$-genus for
$S^1 \backslash M$. If there were a Dirac operator on $S^1 \backslash M$ then
this $\widehat{A}$-genus should be its index. Although we will not actually
construct a Dirac operator on $S^1 \backslash M$, it is nevertheless 
worth considering the topological conditions to have such an operator.
Suppose that $M$ is spin, with a free $S^1$-action.
It does not follow that $S^1 \backslash M$
is spin.  For example, if $M = S^{4k+1}$ has the Hopf action then
$S^1 \backslash M = \C P^{2k}$, which is not spin.
The problem in this case is that the $S^1$-action on the oriented orthonormal
frame bundle of M does not lift to an
$S^1$-action on the principal spin bundle. 
Recall that an $S^1$-action is said to be even if it lifts to
the principal spin bundle and odd if it does not
\cite[p. 295]{Lawson-Michelsohn (1989)}.
We will consider the two cases separately.

\begin{lemma}
Let $M$ be a spin manifold with a fixed spin structure and a semifree
$S^1$-action.
If $F$ is a connected component of
$M^{S^1}$, let $\codim(F)$ be its codimension in $M$. \\
1. If the $S^1$-action is even then $\codim (F) = 2$ or 
$\codim(F) \equiv 0 \mod{4}$.\\
2. If the $S^1$-action is odd then $\codim(F) \equiv 2 \mod{4}$.
\end{lemma}
\begin{pf}
Let $NF$ be the normal bundle to $F$ and let $SNF$ be its sphere bundle, with
fiber $S^{2N+1}$. Then $\codim(F) = 2N + 2$.
As $S^1$ acts trivially on $F$, if the $S^1$-action on $M$
is even (odd) then the Hopf action on $S^{2N+1}$ is even (odd).
(Note that $S^{2N+1}$ has a unique spin structure if $N > 0$.)
If the Hopf action on $S^{2N+1}$ is even then either $N = 0$ and the
spin structure on $S^1$ is the one which does not extend to $D^2$, or
$N$ is odd.
Thus $F$ satisfies conclusion 1. of the 
lemma.
If the Hopf action on $S^{2N+1}$ is odd then $N$ is even, so 
$F$ satisfies conclusion 2. of the lemma.
\end{pf}

\subsubsection{Even Semifree $S^1$-Actions}

Suppose that the spin manifold
$M$ has an even effective $S^1$-action.
Let $SM$ be the spinor bundle of $M$. If $\dim(M) =
4k+1$ then $\dim_{\C} SM = 2^{2k}$. 
If the $S^1$-action is free
then $S^1 \backslash M$ acquires a spin structure, with spinor bundle
$S(S^1 \backslash M) = S^1 \backslash SM$. If the $S^1$-action is semifree,
let $M^{S^1}_{(2)}$ denote the submanifold of $M^{S^1}$ which has
codimension $2$ in $M$. As $M^{S^1}_{(2)}$ appears as a boundary component
in a compactification of $(S^1 \backslash M) - M^{S^1}$, it acquires a
spin structure. Let $D_{M^{S^1}_{(2)}}$ denote the Dirac operator on
$M^{S^1}_{(2)}$.

\begin{definition}
\begin{equation} \label{eqn35}
\hA_{S^1}(M) = \int_{S^1 \backslash M} \hA
\left( T \left( S^1 \backslash M \right) \right) 
+ \frac{1}{2}
\left[ \eta \left( D_{M^{S^1}_{(2)}} \right) + 
\dim \left( \Ker \left( D_{M^{S^1}_{(2)}}
\right) \right) \right].
\end{equation}
\end{definition}

\begin{proposition} \label{prop10}
The number $\hA_{S^1}(M)$ is an integer. 
If $\{g(\epsilon) \}_{\epsilon \in [0,1]}$ is a smooth $1$-parameter family
of $S^1$-invariant metrics on $M$ and $\dim \left( \Ker \left( 
D_{M^{S^1}_{(2)}}
\right) \right)$ is constant in $\epsilon$ then $\hA_{S^1}(M)$ is
constant in $\epsilon$.  
\end{proposition}
\begin{pf}
For small $r > 0$, let $N_r \left( M^{S^1} \right)$ be the $r$-neighborhood
of $M^{S^1}$ in $S^1 \backslash M$. The manifold-with-boundary
$(S^1 \backslash M) - N_r \left( M^{S^1} \right)$ is spin and one can
talk about the index $\Ind_r \in \Z$ of its Dirac operator. By the 
method of proof of Theorem \ref{prop7}, one finds that in $\R / \Z$,
\begin{equation}
\int_{S^1 \backslash M} \hA
\left( T \left( S^1 \backslash M \right) \right) 
+ \lim_{r \rightarrow 0} \frac{1}{2}
\left[ \eta \left( \partial N_r(M^{S^1}) \right) + 
\dim \left( \Ker \left( \partial N_r(M^{S^1})
\right) \right) \right] \equiv 0.
\end{equation}
(The spectral invariants in the above equation are with respect to
Dirac operators.)
Let $F$ be a connected component of $M^{S^1}$ whose codimension in $M$ is 
divisible by four. Then $\partial N_r(F)$ is a fiber bundle whose fiber
is $\C P^{N}$ for some odd $N$. As $\C P^N$ is a spin manifold with positive
scalar curvature, it follows from \cite{Bismut-Cheeger (1989)} that
\begin{equation}
\lim_{r \rightarrow 0} 
\eta \left( \partial N_r(M^{S^1}) \right) = 
\int_F \widetilde{\eta} \wedge \widehat{A}(TF)
\end{equation}
and
\begin{equation}
\lim_{r \rightarrow 0}
\dim \left( \Ker \left( \partial N_r(M^{S^1})
\right) \right) = 0.
\end{equation}
As in the proof of Theorem \ref{prop7}, $\widetilde{\eta} = 0$.

If $F$ is a connected component of $M^{S^1}$ whose codimension in $M$ is two
then $\partial N_r (F)$ is a Riemannian manifold which is topologically the 
same as $F$ and which approaches $F$ metrically as $r \rightarrow 0$.
Thus in $\R / \Z$,
\begin{equation}
\lim_{r \rightarrow 0} \frac{1}{2}
\left[ \eta \left( \partial N_r(F) \right) + 
\dim \left( \Ker \left( \partial N_r(F)
\right) \right) \right] \equiv
\frac{1}{2}
\left[ \eta \left( D_F \right) + 
\dim \left( \Ker \left( D_F
\right) \right) \right].
\end{equation}
It follows that $\widehat{A}_{S^1}(M)$ is an integer.

Let $\{g(\epsilon)\}_{\epsilon \in [0,1]}$ be a family of metrics as in the
statement of the proposition.
Let $I_{S^1}(M)$ denote the first term in the right-hand-side of (\ref{eqn35}).
We first compute 
$I_{S^1}(M) \big|_{\epsilon = 1} - I_{S^1}(M) \big|_{\epsilon = 0}$. 
As in the proof of Theorem 
\ref{prop7}, we compactify $\left( S^1 \backslash M \right) - M^{S^1}$ by
$\cup_F \left( S^1 \backslash SNF \right)$, 
where $F$ ranges over the connected components
of $M^{S^1}$. 
Let $\widehat{\omega}(\epsilon)$ be the connection on
$S^1 \backslash SNF$, as in (\ref{eqn18}).
We can compute $I_{S^1}(M) \big|_{\epsilon = 1} - 
I_{S^1}(M) \big|_{\epsilon = 0}$ as the integral over $\cup_F \left(
S^1 \backslash SNF \right)$ of a transgressed characteristic class.  Namely,
\begin{equation} \label{eqn37}
I_{S^1}(M) \big|_{\epsilon = 1} - I_{S^1}(M) \big|_{\epsilon = 0} =
- \sum_F \int_0^1 \int_{S^1 \backslash SNF} 
\hA \left( \widehat{\Omega}(\epsilon) + d\epsilon \wedge 
\partial_\epsilon \widehat{\omega}  \right). 
\end{equation}
(The minus sign on the right-hand-side of (\ref{eqn37}) comes from the
different orientations of $S^1 \backslash SNF$.)
Let $\widehat{\omega}_V$ be the $i$ and $r$-components of (\ref{eqn18}).
Then from the structure of (\ref{eqn18}),
\begin{align} \label{eqn38}
I_{S^1}(M) \big|_{\epsilon = 1} - I_{S^1}(M) \big|_{\epsilon = 0} =
- \sum_F \int_0^1 \int_F &
\hA \left( R^{TF}(\epsilon) + d\epsilon \wedge 
\partial_\epsilon \omega^{TF}  \right) \wedge \\
& \int_Z 
\hA \left( \widehat{\Omega}_V(\epsilon) + d\epsilon \wedge 
\partial_\epsilon \widehat{\omega}_V  \right). \notag 
\end{align}
Let us write
\begin{align} \label{eqn39}
\hA \left( R^{TF}(\epsilon) + d\epsilon \wedge 
\partial_\epsilon \omega^{TF}  \right) & = 
a_1 + d\epsilon \wedge a_2, \\
\hA \left( \widehat{\Omega}_V(\epsilon) + d\epsilon \wedge 
\partial_\epsilon \widehat{\omega}_V  \right) 
& = b_1 + d\epsilon \wedge b_2, \notag
\end{align}
where $a_1, a_2, b_1$ and $b_2$ depend on $\epsilon$.
Then
\begin{equation} \label{eqn40}
I_{S^1}(M) \big|_{\epsilon = 1} - I_{S^1}(M) \big|_{\epsilon = 0} =
- \sum_F \int_0^1 d\epsilon \wedge
\left( \int_F a_1 \wedge \int_Z b_2 + \int_F a_2 \wedge \int_Z b_1 \right).
\end{equation}
Now $b_1$ and $b_2$ can be computed by equivariant means, and the result
will be a polynomial in the curvature of the principal bundle underlying
$S^1 \backslash SNF$. In particular, they will be even forms.  However,
by parity considerations, $b_2$ is an odd form. Thus $b_2 = 0$ and
\begin{equation} \label{eqn41}
I_{S^1}(M) \big|_{\epsilon = 1} - I_{S^1}(M) \big|_{\epsilon = 0} =
- \sum_F \int_0^1 d\epsilon \wedge
\int_F a_2 \wedge \int_Z b_1.
\end{equation}

From (\ref{eqn39}),
\begin{equation} \label{eqn42}
b_1 = \hA \left( \widehat{\Omega}_V(\epsilon) \right). 
\end{equation}
The Atiyah-Singer families index theorem gives 
an equality in $\HH^{even}(F; \R)$:
\begin{equation} \label{eqn43}
\ch (\Ind (D_Z)) = \int_Z b_1, 
\end{equation}
where $D_Z$ is the family of vertical Dirac operators on the fiber bundle
$S^1 \backslash SNF \rightarrow F$.
If $\dim(Z) > 0$ then $Z$ is a spin manifold with positive scalar curvature
and so 
$\Ind(D_Z) = 0$.
If $\dim(Z) = 0$ then $Z$ is a point and
$\int_Z b_1 = 1$.
Thus
\begin{equation} \label{eqn45}
I_{S^1}(M) \big|_{\epsilon = 1} - I_{S^1}(M) \big|_{\epsilon = 0} =
- \int_0^1 d\epsilon \wedge
\int_{M^{S^1}_{(2)}} a_2.
\end{equation} 

On the other hand, from \cite{Atiyah-Patodi-Singer II (1975)}, as
there is no spectral flow,
\begin{align} \label{eqn46}
\frac{1}{2}
\left[ \eta \left( D_{M^{S^1}_{(2)}} \right) + 
\dim \left( \Ker \left( D_{M^{S^1}_{(2)}}
\right) \right) \right] \Big|_{\epsilon = 1} & -
\frac{1}{2}
\left[ \eta \left(D_{M^{S^1}_{(2)}} \right) + 
\dim \left( \Ker \left( D_{M^{S^1}_{(2)}}
\right) \right) \right] \Big|_{\epsilon = 0} = \\
& \int_0^1 d\epsilon \wedge
\int_{M^{S^1}_{(2)}} a_2. \notag
\end{align} 
The proposition follows.
\end{pf}

\begin{theorem} \label{pscthm}
If $M$ admits an
$S^1$-invariant metric of positive scalar curvature and 
$M^{S^1}_{(2)} = \emptyset$ then $\widehat{A}_{S^1}(M) = 0$.
\end{theorem}
\begin{pf}
We may assume that $M$ has dimension $4k+1$. 
Suppose that it has an $S^1$-invariant Riemannian metric of
positive scalar curvature.
Let $g^*$ be the quotient Riemannian metric on $(S^1 \backslash M) - M^{S^1}$.
Let $l \in C^\infty \left( (S^1 \backslash M) - M^{S^1} \right)$ 
be the function which assigns to a point $x \in (S^1 \backslash M) - M^{S^1}$
the length of the $S^1$-orbit over $x$. Put
$\widetilde{g} = l^{\frac{2}{4k-1}} \: g^*$. Then
$\widetilde{g}$ has positive scalar curvature
\cite[p. 22]{Berard-Bergery (1983)}.

Let us first suppose that for each connected component $F$ of
$M^{S^1}$, a neighborhood of $F$ in $M$ is $S^1$-isometric to $DNF$.
Then as in (\ref{eqn17}), a local orthonormal basis of $1$-forms on
$S^1 \backslash DNF$ for $\widetilde{g}$ is given by
\begin{align} \label{eqn46.1}
\widehat{\tau}^r & = r^{\frac{1}{4k-1}} \: dr, \\
\widehat{\tau}^i & = r^{\frac{4k}{4k-1}} \: \tau^i, \notag \\
\widehat{\tau}^\alpha & = r^{\frac{1}{4k-1}}\tau^\alpha. \notag
\end{align}
Changing variable to $u \: = \: r^{\frac{4k}{4k-1}}$, we obtain the 
local orthonormal
basis
\begin{align} \label{eqn46.2}
\widehat{\tau}^u & = \frac{4k-1}{4k} \: du, \\
\widehat{\tau}^i & = u \: \tau^i, \notag \\
\widehat{\tau}^\alpha & = u^{\frac{1}{4k}}\tau^\alpha. \notag
\end{align} 

Let us consider a more general class of bases given by
\begin{align} \label{eqn46.3}
\widehat{\tau}^u & = \frac{1}{f(u)} \: du, \\
\widehat{\tau}^i & = u \: \tau^i, \notag \\
\widehat{\tau}^\alpha & = u^{\frac{1}{4k}} \tau^\alpha \notag
\end{align} 
for some positive function $f$. Let
$\phi \in C^\infty(0, \infty)$ be a nondecreasing function such that
$\phi(x) = x$ if $x \in (0, \frac12)$ and $\phi(x) = 1$ if $x \ge 1$.  
Given a small $\epsilon > 0$, define 
\begin{equation} \label{eqn46.4}
f(u) \: = \: \frac{4k}{4k-1} \: \phi \left( \frac{u-\epsilon}{\epsilon^{1/2}}
\right)
\end{equation}
for $u > \epsilon$.
Then one can check that the metric for which (\ref{eqn46.3}) 
is an orthonormal basis is complete with  positive scalar
curvature. In effect, a change of variable to $s \: = \: - \: 
\ln(u- \epsilon)$ shows that the
metric is asymptotically cylindrical, with cross-section 
$S^1 \backslash SNF$ having $\C P^N$-fibers of diameter proportionate to 
$\epsilon$ and base $F$ of diameter proportionate to
$\epsilon^{\frac{1}{4k}}$. As $N > 0$, the positive scalar curvature of the
$\C P^N$ fibers  
ensures that the metric will have positive scalar curvature for
small $\epsilon$.  Truncate the cylinder at a large distance and
smooth the metric to a product near the boundary, while keeping positive scalar
curvature.  Let $N_\epsilon$ denote the corresponding manifold-with-boundary.
Applying the Atiyah-Patodi-Singer theorem and the Lichnerowicz vanishing
theorem to $N_\epsilon$, we obtain
\begin{equation} \label{eqn46.5}
0 \: = \: \int_{N_\epsilon} \hA
\left( T N_\epsilon \right) 
+ \frac{1}{2}
\: \eta \left( D_{\partial N_\epsilon} \right).
\end{equation}
Now 
\begin{equation} \label{eqn46.6}
\lim_{\epsilon \rightarrow 0} \int_{N_\epsilon} \hA
\left( T N_\epsilon \right)  \: = \: \int_{S^1 \backslash M} \hA
\left( T (S^1 \backslash M \right). 
\end{equation}
As in the proof of Proposition \ref{prop10}, since
$M^{S^1}_{(2)} = \emptyset$,
\begin{equation} \label{eqn46.7}
\lim_{\epsilon \rightarrow 0} \eta \left( D_{\partial N_\epsilon} \right) = 0.
\end{equation}
The proposition follows in this case.

In general, a neighborhood of $F$ in $M$ may not be $S^1$-isometric to $DNF$.
Nevertheless, we can use the distance function from $F$ to write 
$\widetilde{g}$ as $\left( \frac{4k-1}{4k} \right)^2 \: du^2 + h(u)$ where 
for $u > 0$, $h(u)$ is a metric on $S^1 \backslash SNF$. For small $u$, 
$\widetilde{g}$ will be well-approximated by the metric of the form 
(\ref{eqn46.2}). 
Then
we can deform $\widetilde{g}$ for small $u$ to obtain a metric of positive
scalar curvature and precisely
of the form (\ref{eqn46.2}) for small $u$, to which we can
apply the previous argument. 
\end{pf}
{\bf Remark :} Suppose that $M^{S^1}$ has codimension two in $M$.
Then the orthonormal frame (\ref{eqn46.2}) becomes
\begin{align} \label{eqn46.8}
\widehat{\tau}^u & = \frac{4k}{4k-1} \: du, \\
\widehat{\tau}^\alpha & = u^{\frac{1}{4k}}\tau^\alpha. \notag
\end{align} 
We no longer have the benefit of the positive scalar curvature 
coming from $\C P^N$.
Metrically with respect to $\widetilde{g}$,  $S^1 \backslash M$ has a 
``puffy'' cone over $M^{S^1}_{(2)}$. If one could prove an index theorem
for Dirac operators on such spaces, along with a vanishing theorem in the
case of positive scalar curvature, one could remove the codimension
restriction in Theorem \ref{pscthm}.

\subsubsection{Even or Odd Semifree $S^1$-Actions}

Suppose that the spin manifold $M$ has an $S^1$-action which is even or odd.
If the $S^1$-action is free then $S^1 \backslash M$ may not have a spin
structure, but it always has a canonical $\spin^c$ structure.
Namely, if the $S^1$-action is even, put $S(S^1 \backslash M) =
\C \times_{S^1} SM$, where $\C$ has the standard $S^1$-action.
If the $S^1$-action is odd, let $\widehat{S^1}$ be the double cover of
$S^1$. It acts on $M$ through the quotient map 
$\widehat{S^1} \rightarrow S^1$. Consider the standard action of 
$\widehat{S^1}$ on $\C$. 
The infinitesimal action of $u(1)$ on $\C \times SM$ integrates to an
$\widehat{S^1}$-action, so we can put 
$S(S^1 \backslash M) =
\C \times_{\widehat{S^1}} SM$. In either case, $S(S^1 \backslash M)$ is
the spinor bundle of a $\spin^c$ structure on $S^1 \backslash M$.

Now suppose that the $S^1$-action is effective and semifree.
\begin{lemma} 
$M^{S^1}$ is $\spin^c$.
\end{lemma}
\begin{pf}
Let $F$ be a connected component of $M^{S^1}$, with normal bundle $NF$.
We know that $F$ is oriented. As the total space
$NF$ is diffeomorphic to a neighborhood of $F$ in $M$, $TNF$ inherits a
spin structure.  Let $p : NF \rightarrow F$ be projection to the base.
Then $TNF = p^*NF \oplus p^* TF$.
As $NF$ has a complex structure, it has a canonical $\spin^c$ structure.
Then $p^* TF$ acquires a $\spin^c$ structure, and so does $TF$.
\end{pf}

Let $\xi^{S^1 \backslash M}$ 
be the complex line bundle on $\left( S^1 \backslash M \right) -
M^{S^1}$ associated to the $\spin^c$ structure. It has an induced connection
$\nabla^{\xi^{S^1 \backslash M}}$. Let $c_1 \left( \xi^{S^1 \backslash M}
\right) \in \Omega^2 \left( (S^1 \backslash M) - M^{S^1} \right)$ be the 
corresponding characteristic form. Let $D_{M^{S^1}}$ be the $\spin^c$ Dirac
operator on $M^{S^1}$. 

\begin{definition}
\begin{equation} \label{eqn35.}
\hcA_{S^1}(M) = \int_{S^1 \backslash M} \hA
\left( T \left( S^1 \backslash M \right) \right) 
\wedge e^{\frac{c_1({\xi^{S^1 \backslash M}})}{2}} + \frac{1}{2}
\left[ \eta \left( D_{M^{S^1}} \right) + \dim \left( \Ker \left( D_{M^{S^1}}
\right) \right) \right].
\end{equation}
\end{definition}

\begin{proposition} \label{prop11}
The number $\hcA_{S^1}(M)$ is an integer.
If $\{g(\epsilon) \}_{\epsilon \in [0,1]}$ is a smooth $1$-parameter family
of $S^1$-invariant metrics on $M$ and $\dim \left( \Ker \left( D_{M^{S^1}}
\right) \right)$ is constant in $\epsilon$ then $\hcA_{S^1}(M)$ is
constant in $\epsilon$. 
\end{proposition}
\begin{pf}
For notational convenience, put
\begin{equation}
\hcA(TX) = \hA(TX) \wedge e^{\frac{c_1({\xi^X})}{2}}.
\end{equation}
For small $r > 0$, let $N_r \left( M^{S^1} \right)$ be the $r$-neighborhood
of $M^{S^1}$ in $S^1 \backslash M$. The manifold-with-boundary
$(S^1 \backslash M) - N_r \left( M^{S^1} \right)$ is $\spin^c$ and one can
talk about the index $\Ind_r \in \Z$ of its Dirac operator. By the 
method of proof of Theorem \ref{prop7}, one finds that in $\R / \Z$,
\begin{equation}
\int_{S^1 \backslash M} \hcA
\left( T \left( S^1 \backslash M \right) \right) 
+ \lim_{r \rightarrow 0} \frac{1}{2}
\left[ \eta \left( \partial N_r(M^{S^1}) \right) + 
\dim \left( \Ker \left( \partial N_r(M^{S^1})
\right) \right) \right] \equiv 0.
\end{equation}
(The spectral invariants in the above equation are with respect to
$\spin^c$ Dirac operators.)
Let $F$ be a connected component of $M^{S^1}$.
Then $\partial N_r(F)$ is a fiber bundle over $F$.
In terms of the complex structure on 
a fiber $Z$, we can write 
\begin{equation} \label{eqn44.}
D_Z = \overline{\partial} + \overline{\partial}^* : \Omega^{0, even}(Z)
\rightarrow \Omega^{0, odd}(Z).
\end{equation}
As the fiber $Z$ is a complex projective space,  
$\Ker(D_Z) = \Omega^{0,0}(Z) = \C$ consists of the constant functions
on the fibers and $\Ker(D_Z^*) = 0$.  
Hence $\Ind(D_Z)$ is a trivial complex line bundle on $F$.
It follows from \cite{Dai (1991)} that in $\R / \Z$,
\begin{equation}
\lim_{r \rightarrow 0} 
\frac{1}{2}
\left[ \eta \left( \partial N_r(F) \right) + 
\dim \left( \Ker \left( \partial N_r(F)
\right) \right) \right] \equiv
\frac{1}{2}
\left[ \int_F \widetilde{\eta} \wedge \hcA(TF) + \eta(D_F)
 + \dim \left( \Ker \left( D_F
\right) \right) \right].
\end{equation}
As in the proof of Theorem \ref{prop7}, $\widetilde{\eta} = 0$.
Thus 
$\hcA_{S^1}(M)$
is an integer.

The rest of the proof of Proposition \ref{prop11} is similar to that of
Proposition \ref{prop10}. We omit the details.
\end{pf}

\subsection{General $S^1$-Actions}
\label{General $S^1$-Actions}

Let $S^1$ act effectively on $M$. There are suborbifolds ${\cal O}$ of
$(S^1 \backslash M) - M^{S^1}$ defined as in Subsection 
\ref{Fixed-point-free Actions}. 
\begin{proposition} \label{prop12}
\begin{equation} \label{eqn47}
\sigma_{S^1}(M) = \sum_{\cal O} \frac{1}{m_{\cal O}}
\int_{\cal O} L({\cal O}) + \eta \left(M^{S^1} \right).
\end{equation}
\end{proposition}
\begin{pf}
The proof is a combination of those of 
(\ref{eqn12}) and Theorem \ref{prop7}. Let $F$ be a
connected component of $M^{S^1}$. Let $NF$ be the normal bundle of $F$ in
$M$. It has an $S^1$-action by orthogonal automorphisms, which is 
fixed-point-free on $NF - F$. Let $SNF$ be the sphere bundle of $NF$.
Then $S^1 \backslash SNF$ is an orbifold. For $r > 0$, let $N_r(F)$ be
the $r$-neighborhood of $F$ in $S^1 \backslash M$. Then for small $r$,
$\partial N_r(F)$ is an orbifold. We define the $\eta$-invariant of
$\partial N_r(F)$ using the tangential signature operator on 
orbifold-differential forms on $\partial N_r(F)$, i.e. 
$S^1$-basic differential forms on the preimage of $\partial N_r(F)$ in $M$.
Then the method of proof of Theorem \ref{prop7} goes through with minor
changes.
\end{pf}

\section{Equivariant Higher Indices} \label{Equivariant Higher Indices}
\subsection{Equivariant Novikov Conjectures}
\label{Equivariant Novikov Conjectures}
Let $M^n$ be a closed oriented connected
manifold. Let $\Gamma^\prime$ be a countable discrete group and
let $\rho : \pi_1(M) \rightarrow \Gamma^\prime$ be a surjective homomorphism.
There is an induced continuous
map $\nu : M \rightarrow B\Gamma^\prime$, defined up to
homotopy. Let $L \in \HH_{n-4*}(M; \Q)$ 
be the homology $L$-class of $M$, i.e.
the Poincar\'e dual of the cohomology $L$-class. The
Novikov Conjecture hypothesizes that 
$\nu_*(L) \in \HH_{n-4*}(B\Gamma^\prime; \Q)$ is
an oriented homotopy invariant of $M$. Another way to state this is to
let $D \in \KO_n(M)$ be the $\KO$-homology class of the signature operator.
The Novikov Conjecture
says that $\nu_*(D) \otimes_\Z 1 \in \KO_n(B\Gamma^\prime) \otimes_\Z \Q$ 
should be an oriented homotopy invariant of $M$. 
It is usually assumed that $\Gamma^\prime = \pi_1(M)$, although this is 
not necessary.

Now suppose that a compact Lie group $G$ acts on $M$ in an 
orientation-preserving way. One would like to extend the Novikov
Conjecture to the $G$-equivariant setting.  One approach is to extend the
classifying space construction. The idea is that $B\pi_1(M)$ has exactly the 
information about $\pi_0(M)$ and $\pi_1(M)$. In the equivariant case one wants
a space with a $G$-action, constructed from the data 
$\{\pi_0(M^H)\}$ and $\{\pi_1(M^H)\}$ as $H$ runs over the closed subgroups
of $G$. Such a space $B\pi(M)$ is constructed in \cite{May (1990)}. It
has the property that each connected component of $B\pi(M)^H$ is 
aspherical, and
there is a $G$-map $\nu : M \rightarrow B\pi(M)$, unique up to $G$-homotopy, 
which induces an isomorphism from $\pi_0(M^H)$ to $\pi_0(B\pi(M)^H)$
and an isomorphism on $\pi_1$ of each
connected component of $M^H$. Choosing a $G$-invariant Riemannian
metric on $M$, there is a $G$-invariant signature operator $D \in 
\KO_n^G(M)$. Then one Equivariant Novikov Conjecture would be that
$\nu_*(D) \otimes_\Z 1  \in \KO_n^G(B\pi(M)) \otimes_\Z \Q$ is an oriented 
$G$-homotopy invariant of $M$ \cite{Rosenberg-Weinberger (1990)}.

As was pointed out in \cite[p. 31]{Rosenberg-Weinberger (1990)}, this 
conjecture is false in the case of free $S^1$-actions. In that case
$B\pi(M) = S^\infty$, $\KO_n^G(B\pi(M)) = \KO_{n-1} 
\left( \C P^\infty \right)$ 
and $\KO_n^G(B\pi(M)) \otimes_\Z \Q = 
\HH_{n-1-4*} \left( \C P^\infty ; \Q \right)$.
The principal $S^1$-bundle $M$ is classified by a map $f : (S^1 \backslash M)
\rightarrow \C P^\infty$, 
and $\nu_*(D) \otimes_\Z 1 = f_*(L(S^1 \backslash M)) \in
\HH_{n-1-4*} \left( \C P^\infty ; \Q \right)$.
If $X$ is a homotopy-$\C P^N$, let $M$ be the total space of the $S^1$-bundle
associated to the standard generator of $\HH^2 \left( X; \Z \right) =
\HH^2 \left( \C P^N; \Z \right)$.
Then $\nu_*(D) \otimes_\Z 1$ can be identified with 
the rational homology $L$-class of $X$.
If $N > 2$ then it follows from surgery theory that
there is an infinite number of nonhomeomorphic 
homotopy-$\C P^N$'s $\{X_i\}_{i=1}^\infty$ with distinct
rational homology $L$-classes. The $S^1$-actions on the
corresponding homotopy-spheres $\{M_i\}_{i=1}^\infty$ will be mutually
homotopy equivalent, showing the falsity of the conjecture.
The rest of \cite{Rosenberg-Weinberger (1990)} is devoted to looking at the
conjecture under some finiteness assumptions on $B\pi(M)$.

Another Equivariant Novikov Conjecture uses the classifying space
$\underline{E} G^\prime$ for proper $G^\prime$-actions, where 
$G^\prime$ is a Lie group with a countable number of connected components
\cite{Baum-Connes-Higson (1994)}.
Let $\Gamma^\prime$ and $\rho$ be as above.
There is an induced connected normal $\Gamma^\prime$-covering 
$M^\prime$ of $M$.
Let $\pi : M^\prime \rightarrow M$ be the projection map.
Define a group $G^\prime$ by
\begin{equation} \label{eqn48}
G^\prime = \{ (\phi, g) \in \Diff \left( M^\prime \right) \times G :
\pi \circ \phi = g \cdot \pi \}. 
\end{equation}
There is a $G^\prime$-invariant signature operator $D \in 
\KO_n^{G^\prime}(M^\prime)$. The conjecture states that
$\nu_*(D) \in \KO_n^{G^\prime}(\underline{E} G^\prime)$ 
is an oriented $G$-homotopy invariant of $M$
\cite{Baum-Connes-Higson (1994)},
\cite[Proposition 2.10]{Rosenberg-Weinberger (1990)}. 

This conjecture is very reasonable. However, it seems to 
be more useful when $G$ is finite.  
Suppose, for example, that $G = S^1$ and $\Gamma^\prime =
\{e\}$. Then $G^\prime = S^1$, $\underline{E}S^1$ is a point and 
if $n$ is divisible by four then $\KO_n^{S^1}(\pt.)$ is a countable sum of 
$\Z$'s, while it vanishes rationally otherwise.
If $n$ is divisible by four then the only information in $\nu_n(D)
\in \KO_n^{S^1}(\pt.)$ is the ordinary signature
of $M$. If $S^1$ acts freely on $M$ then $M = \partial \left( D^2 \times_{S^1}
M \right)$ and so its signature vanishes.  Thus in the
case of free $S^1$-actions, the second Equivariant Novikov Conjecture is 
true but vacuous.

In order to construct higher signatures of $S^1 \backslash M$, we will
use the higher eta-invariant of \cite{Lott (1992a)}. We now recall the
construction of \cite{Lott (1992a)}, with some modifications.
We will let groups act on the left, as in \cite{Lott (1996)}, instead of
on the right, as in \cite{Lott (1992a)}. The differential form conventions
will be as in \cite{Lott (1996)}.

\subsection{Higher Eta-Invariant}
\label{Higher Eta-Invariant}

Let $\Gamma$ be a finitely generated discrete group and let $C^*_r \Gamma$
be the reduced group $C^*$-algebra. 
\begin{assumption} \label{ass1}
There is a Fr\'echet locally $m$-convex algebra ${\frak B}$ such
that\\
1. $\C \Gamma \subset {\frak B} \subset C^*_r \Gamma$.\\
2. ${\frak B}$ is stable under the holomorphic functional calculus
in $C^*_r \Gamma$.\\
3. For each $\tau \in \HH^q(\Gamma; \C)$, there is a representative cocycle
$\tau \in Z^q(\Gamma; \C)$ such that the ensuing cyclic cocycle $Z_\tau
\in ZC^q(\C \Gamma)$ extends to a continuous cyclic cocycle on 
${\frak B}$.
\end{assumption}

It is know that such ``smooth subalgebras'' ${\frak B}$ exist if
$\Gamma$ is virtually nilpotent or Gromov-hyperbolic 
\cite[Section III.5]{Connes (1994)}, \cite{Ji (1992)}.

Let $F$ be a closed oriented Riemannian manifold of dimension $n$. Let 
$\rho : \pi_1(F) \rightarrow
\Gamma$ be a surjective homomorphism. There is an induced connected normal
$\Gamma$-covering $F^\prime$ of $F$, on which $g \in \Gamma$ acts on
the left by $L_g \in \Diff \left(F^\prime \right)$. Let 
$\pi : F^\prime \rightarrow F$ be the projection map.  

Put ${\cal D} = {\frak B} \times_\Gamma F^\prime$, a 
${\frak B}$-vector bundle on $F$, and put $\overline{{\cal D}} =
(C^*_r \Gamma) \times_\Gamma F^\prime$, a $C^*_r \Gamma$-vector bundle on $F$.
Both ${\cal D}$ and $\overline{\cal D}$ are local systems.

\begin{assumption} \label{ass2}
If $n$ is even then $\HH^{\frac{n}{2}}(F; \overline{\cal D}) = 0$. 
If $n$ is odd then 
$\HH^{\frac{n \pm 1}{2}}(F; \overline{\cal D}) = 0$.
\end{assumption}

The cohomology involved in Assumption \ref{ass2} 
is ordinary unreduced cohomology; 
that is, we quotient by $\Image(d)$, not its closure. Equivalent formulations
are : \\
1. If $n$ is even then the spectrum of the $L^2$-Laplacian on $F^\prime$
is strictly positive in degree $\frac{n}{2}$. If $n$ is odd then
the spectrum of the $L^2$-Laplacian on $F^\prime$
is strictly positive in degrees $\frac{n \pm 1}{2}$. \\
2. If $n$ is even then the Laplacian on
$\Omega^{\frac{n}{2}}(F; \overline{\cal D})$ is invertible. If $n$ is odd 
then the Laplacians on
$\Omega^{\frac{n \pm 1}{2}}(F; \overline{\cal D})$ are invertible.\\
3. If $n$ is even then the Laplacian on
$\Omega^{\frac{n}{2}}(F; {\cal D})$ is invertible. If $n$ is odd 
then the Laplacians on
$\Omega^{\frac{n \pm 1}{2}}(F; {\cal D})$ are invertible.\\
4. If $n$ is even then $\HH^{\frac{n}{2}}(F; {\cal D}) = 0$. If $n$ is odd
then $\HH^{\frac{n \pm 1}{2}}(F; {\cal D}) = 0$. \\

We use the notions of Hermitian complex and regular Hermitian complex
from \cite{Kaminker-Miller (1985)} and \cite{Lusztig (1971)}.
Using \cite[Section 4.1 and Proposition 10]{Lott (1996)}, one can
generalize the results of \cite{Kaminker-Miller (1985)} from 
$C^*_r \Gamma$-complexes to ${\frak B}$-complexes.
 
\begin{proposition} \label{prop13}
There is a cochain complex $W^*$ of finitely-generated projective
${\frak B}$-modules such that\\
1. $W^*$ is a regular Hermitian complex.\\
2. $W^{\frac{n}{2}} = 0$ if $n$ is even and
$W^{\frac{n \pm 1}{2}} = 0$ if $n$ is odd.\\
3. The complex $\Omega^*(F; {\cal D})$ of smooth ${\cal D}$-valued differential
forms on $F$ is homotopy equivalent to $W^*$.
\end{proposition}
\begin{pf}
We will implicitly use results from  \cite[Proposition 10 and
Section 6.1]{Lott (1996)} concerning spectral analysis involving ${\frak B}$.
First, let $K$ be a triangulation of $F$. Then
$\Omega^*(F; {\cal D})$ is homotopy equivalent to the simplicial cochain
complex $C^*(K; {\cal D})$. The latter is a Hermitian complex of 
finitely-generated free ${\frak B}$-modules.  By
\cite[Proposition 2.4]{Kaminker-Miller (1985)}, 
it is homotopy equivalent to a regular Hermitian complex $V^*$ of
finitely-generated projective ${\frak B}$-modules. 
Suppose that $n$ is even. We have $\HH^{\frac{n}{2}}(V) = 0$.
Put 
\begin{equation} \label{eqn49}
W^i =
\begin{cases}
V^i & \text{if $i < \frac{n}{2} - 1$,} \\
\Ker \left(d : V^{\frac{n}{2}-1} \longrightarrow V^{\frac{n}{2}} \right) &
\text{if $i = \frac{n}{2} - 1$,} \\
0 & \text{if $i = \frac{n}{2}$,} \\
\Image \left(d : V^{\frac{n}{2}} \longrightarrow 
V^{\frac{n}{2}+1} \right)^\perp
& \text{if $i = \frac{n}{2} + 1$,} \\
V^i & \text{if $i > \frac{n}{2} + 1$.} \\
\end{cases}
\end{equation}
Then $W^*$ is a regular Hermitian complex.
There are homotopy equivalences $V^* \longrightarrow W^*$ and
$W^* \longrightarrow V^*$
 given by \\
\begin{equation} \label{eqn50}
\begin{array}{cccccccccccc}
\ldots \longrightarrow & 
V^{\frac{n}{2}-2} & \longrightarrow & V^{\frac{n}{2}-1} & \longrightarrow & 
V^{\frac{n}{2}} & \longrightarrow & V^{\frac{n}{2}+1} & \longrightarrow &
V^{\frac{n}{2}+2} & \longrightarrow & \ldots \\ 
&&&&&&&&&&& \\
 & \Id. \downarrow & & p \downarrow & & 0 \downarrow & & 
p \downarrow & & \Id. \downarrow & &
\\ 
&&&&&&&&&&& \\ 
\ldots \longrightarrow & 
W^{\frac{n}{2}-2} & \longrightarrow & W^{\frac{n}{2}-1} & \longrightarrow & 
0 & \longrightarrow & W^{\frac{n}{2}+1} & \longrightarrow &
W^{\frac{n}{2}+2} & \longrightarrow & \ldots 
\end{array}
\end{equation}
and
\begin{equation} \label{eqn51}
\begin{array}{cccccccccccc}
\ldots \longrightarrow & 
W^{\frac{n}{2}-2} & \longrightarrow & W^{\frac{n}{2}-1} & \longrightarrow & 
0 & \longrightarrow & W^{\frac{n}{2}+1} & \longrightarrow &
W^{\frac{n}{2}+2} & \longrightarrow & \ldots \\
&&&&&&&&&&& \\ 
 & \Id. \downarrow & & i \downarrow & & 0 \downarrow & & 
i \downarrow & & \Id. \downarrow & &
\\
&&&&&&&&&&& \\ 
\ldots \longrightarrow & 
V^{\frac{n}{2}-2} & \longrightarrow & V^{\frac{n}{2}-1} & \longrightarrow & 
V^{\frac{n}{2}} & \longrightarrow & V^{\frac{n}{2}+1} & \longrightarrow &
V^{\frac{n}{2}+2} & \longrightarrow & \ldots, 
\end{array}
\end{equation}
where $p$ denotes orthogonal projection and
$i$ is inclusion. The cochain homotopy operators are
\begin{equation} \label{eqn52}
\ldots \longleftarrow 
V^{\frac{n}{2}-2} \stackrel{0}{\longleftarrow} V^{\frac{n}{2}-1} 
\stackrel{d^* \triangle^{-1}}{\longleftarrow}  
V^{\frac{n}{2}}  \stackrel{\triangle^{-1} d^*}{\longleftarrow}  
V^{\frac{n}{2}+1}
\stackrel{0}{\longleftarrow}
V^{\frac{n}{2}+2}  \longleftarrow \ldots
\end{equation}
and
\begin{equation} \label{eqn53}
\ldots \longleftarrow 
W^{\frac{n}{2}-2} \stackrel{0}{\longleftarrow} W^{\frac{n}{2}-1} 
\stackrel{0}{\longleftarrow}  0  \stackrel{0}{\longleftarrow}  
W^{\frac{n}{2}+1}
\stackrel{0}{\longleftarrow}
W^{\frac{n}{2}+2} \longleftarrow \ldots
\end{equation}

If $n$ is odd, we have $\HH^{\frac{n \pm 1}{2}}(V) = 0$. Put 
\begin{equation} \label{eqn54}
W^i =
\begin{cases}
V^i & \text{if $i < \frac{n-3}{2}$,} \\
\Ker \left(d : V^{\frac{n-3}{2}} \longrightarrow V^{\frac{n-1}{2}} \right) &
\text{if $i = \frac{n-3}{2}$,} \\
0 & \text{if $i = \frac{n \pm 1}{2}$,} \\
\Image \left(d : V^{\frac{n+1}{2}} \longrightarrow 
V^{\frac{n+3}{2}} \right)^\perp
& \text{if $i = \frac{n+3}{2}$,} \\
V^i & \text{if $i > \frac{n+3}{2}$.} \\
\end{cases}
\end{equation}
Then $W^*$ is a regular Hermitian complex.
There are homotopy equivalences $V^* \longrightarrow W^*$ and
$W^* \longrightarrow V^*$ given by \\
\begin{equation} \label{eqn55}
\begin{array}{cccccccccc}
\ldots \longrightarrow & 
V^{\frac{n-3}{2}} & \longrightarrow & V^{\frac{n-1}{2}} & \longrightarrow & 
V^{\frac{n+1}{2}} & \longrightarrow & V^{\frac{n+3}{2}} & 
\longrightarrow & \ldots 
\\ 
&&&&&&&&& \\ 
 & p \downarrow & & 0 \downarrow & & 0 \downarrow & &  
p \downarrow & &
\\ 
&&&&&&&&& \\ 
\ldots \longrightarrow & W^{\frac{n-3}{2}} & \longrightarrow & 0 & 
\longrightarrow & 
0 & \longrightarrow & W^{\frac{n+3}{2}} & \longrightarrow & \ldots 
\end{array}
\end{equation}
and
\begin{equation} \label{eqn56}
\begin{array}{cccccccccc}
\ldots \longrightarrow & W^{\frac{n-3}{2}} & \longrightarrow & 0 & 
\longrightarrow & 
0 & \longrightarrow & W^{\frac{n+3}{2}} & \longrightarrow & \ldots \\
&&&&&&&&& \\ 
 & i \downarrow & & 0 \downarrow & & 0 \downarrow & &  
i \downarrow & & \\ 
&&&&&&&&& \\ 
\ldots \longrightarrow & 
V^{\frac{n-3}{2}} & \longrightarrow & V^{\frac{n-1}{2}} & \longrightarrow & 
V^{\frac{n+1}{2}} & \longrightarrow & V^{\frac{n+3}{2}} & 
\longrightarrow & \ldots 
\end{array}
\end{equation}
The cochain homotopy operators are
\begin{equation} \label{eqn57}
\ldots \stackrel{0}{\longleftarrow} 
V^{\frac{n-3}{2}} \stackrel{{d^*}{\triangle}^{-1}}{\longleftarrow} 
V^{\frac{n-1}{2}} \stackrel{{d^*}{\triangle}^{-1}}{\longleftarrow}  
V^{\frac{n+1}{2}}  \stackrel{{\triangle}^{-1}{d^*}}{\longleftarrow}  
V^{\frac{n+3}{2}} \stackrel{0}{\longleftarrow} \ldots
\end{equation}
and
\begin{equation} \label{eqn58}
\ldots \stackrel{0}{\longleftarrow} 
W^{\frac{n-3}{2}} \stackrel{0}{\longleftarrow} 
0 \stackrel{0}{\longleftarrow}  
0 \stackrel{0}{\longleftarrow}  
W^{\frac{n+3}{2}} \stackrel{0}{\longleftarrow} \ldots
\end{equation}
The proposition follows.
\end{pf}

We briefly review some notation from \cite{Lott (1992a)} and 
\cite{Lott (1996)}. 
Let $\Omega^*({\frak B})$ be the universal graded differential algebra
of ${\frak B}$ and let $\overline{\Omega}^*({\frak B})$ be the quotient
by (the Fr\'echet closure of) the graded commutator. Let
$\overline{\HH}^*({\frak B})$ denote the cohomology of the complex
$\overline{\Omega}^*({\frak B})$.
If $E$ is a complex vector bundle on $F$, put ${\cal E} =
{\cal D} \otimes E$. There is a
bigraded complex 
$\Omega^{*,*} \left( F, {\frak B} \right)$ which, roughly speaking,
consists of differential forms on $F$ along with noncommutative differential 
forms on ${\frak B}$. 

Let $h \in C^\infty_0(F^\prime)$ be a real-valued function satisfying 
$\sum_{g \in \Gamma} L_g^* h = 1$. One obtains a connection
\begin{equation} \label{eqn59}
\nabla^{{\cal D}} : C^\infty \left( F; {\cal D} \right)
\rightarrow \Omega^{1,0} \left( F, {\frak B}; {\cal D} \right)
\oplus \Omega^{0,1} \left( F, {\frak B}; {\cal D} \right)
\end{equation}
on ${\cal D}$. The $(1,0)$-part of the connection comes from the flat
structure of ${\cal D}$ as a vector bundle on $F$. The $(0,1)$-part of
the connection is constructed using $h$. Locally on $F$, using the flat
structure on ${\cal D}$, one can write
\begin{equation} \label{eqn60}
\nabla^{\cal D} = \sum_{\mu = 1}^{dim(F)} dx^\mu \: \partial_\mu +
\sum_{g \in \Gamma} dg \: \nabla_g.
\end{equation}

Suppose that $\dim(F)$ is even. Take $E = \Lambda^* \left( T^* F \right)$,
a vector bundle on $F$ with a $\Z_2$-grading coming from Hodge duality. 
The signature operator
$d + d^* : C^\infty(F; E) \rightarrow C^\infty(F; E)$ couples to 
${\cal D}$ to give a Dirac-type operator
\begin{equation} \label{eqn61}
Q : C^\infty(F; {\cal E}) \rightarrow C^\infty(F; {\cal E})
\end{equation}
which commutes with the left-action of ${\cal B}$.
We can ``quantize'' the $dx^\mu$-variables in
(\ref{eqn60}) to obtain a superconnection
\begin{equation} \label{eqn62}
D : C^\infty \left(F; {\cal E} \right)
\rightarrow C^\infty \left(F; {\cal E} \right) \oplus 
\Omega^{0,1} \left(F, {\frak B}; {\cal E} \right)
\end{equation}
 given by
\begin{equation} \label{eqn63}
D = Q + \nabla^{{\cal E},0,1}.
\end{equation}
Given $s > 0$, we rescale the Clifford variables in (\ref{eqn63}) to obtain
\begin{equation} \label{eqn64}
D_s = s Q + \nabla^{{\cal E},0,1}.
\end{equation}
It extends by Leibniz' rule to an odd map
\begin{equation} \label{eqn65}
D_s : \Omega^{0,*} \left( F, {\frak B}; {\cal E} \right)
\rightarrow \Omega^{0,*} \left(F, {\frak B}; {\cal E} \right).
\end{equation}
(This is like a superconnection on a fiber bundle whose base is the
noncommutative space specified by ${\frak B}$.)
Using the supertrace $\TR_s$ on integral operators on $F$,
one can define
\begin{equation} \label{eqn66}
\TR_s \left( e^{-  D_s^2}
\right) \in \overline{\Omega}^{even} \left( {\frak B} \right).
\end{equation}
A form of the local Atiyah-Singer index theorem says
\begin{equation} \label{eqn67}
\lim_{s \rightarrow 0} \TR_s \left( e^{-  D_s^2} \right) =
\int_F L(TF) \wedge \ch \left( \nabla^{\cal D} \right),
\end{equation}
where $\ch \left( \nabla^{\cal D} \right) \in \Omega^* \left(F; 
\overline{\Omega}^* ({\frak B}) \right)$ is the Chern character.

Consider
\begin{equation} \label{eqn68}
\TR_s \left( Q \: e^{-  D_s^2}
\right) \in \overline{\Omega}^{odd} \left( {\frak B} \right).
\end{equation}
We would like to define the noncommutative eta-form by
\begin{equation} \label{eqn69}
\int_0^\infty \TR_s \left( Q \: e^{-  D_s^2} \right) ds.
\end{equation}
As shown in \cite[Proposition 26]{Lott (1992a)}, there is no problem
with the small-$s$ integration.  In \cite[Section 4.7]{Lott (1992a)} we
argued that the large-$s$ integration is also well-defined, because of
Hodge duality.  However, Eric Leichtnam and Paolo Piazza pointed out to
me that there are technical problems with the argument in 
\cite[Section 4.7]{Lott (1992a)}. Consequently, we do not know whether or
not the integral in (\ref{eqn69}) is convergent for large-$s$. We now present
a way to get around this problem.

Take $n$ of either parity.
For $-1 \le i \le n+1$, put 
\begin{equation} \label{eqn70}
\widehat{W}^i = 
\begin{cases}
W^{i+1} & \text{if $-1 \le i < \frac{n}{2}$,} \\
0 & \text{ if $n$ is even and $i =
\frac{n}{2}$,} \\
W^{i-1}
& \text{if $\frac{n}{2} < i \le n + 1$.} \\
\end{cases}
\end{equation}
and
\begin{equation} \label{eqn71}
C^i = \Omega^i(F; {\cal D}) \oplus \widehat{W}^i.
\end{equation}
Let $f : \Omega^*(F; {\cal D}) \rightarrow W^*$  
be a homotopy equivalence of
Hermitian complexes. 
Let $g : W^* \rightarrow \Omega^*(F; {\cal D})$ be the adjoint of $f$
with respect to the nondegenerate Hermitian form $H_W : W^i \otimes
W^{n-i} \rightarrow {\frak B}$. 
Given $\epsilon \in \R$, define a differential $d^C$ on $C^*$ by
\begin{equation} \label{eqn72}
d^C_i = 
\begin{cases}
\begin{pmatrix}
d & \epsilon g \\
0 & -d
\end{pmatrix}
& \text{if $i < \frac{n}{2}$,} \\
& \\
\: \: d & \text{ if $n$ is even and $i = \frac{n}{2}$,} \\
& \\
\begin{pmatrix}
d & 0 \\
\epsilon f & -d
\end{pmatrix}
& \text{if $i > \frac{n}{2}$.} \\
\end{cases}
\end{equation}

There is a nondegenerate form $H$ on $C^*$ given by
\begin{equation} \label{eqn73}
H((\omega, w), (\omega^\prime, w^\prime)) =
\int_F \omega \wedge \overline{\omega^\prime} + (-1)^{i+1} \: H_W(w, w^\prime)
\end{equation}
for $(\omega, w) \in C^i$, $(\omega^\prime, w^\prime) \in C^{n-i}$,
if $i < \frac{n}{2}$, and
\begin{equation} \label{eqn74}
H(\omega, \omega^\prime) = \int_F \omega \wedge \overline{\omega^\prime}
\end{equation} 
if $\omega, \omega^\prime \in \Omega^{\frac{n}{2}}(F; {\cal D})$.
Then one can check that $C^*$ is a regular Hermitian complex.

If $\epsilon \ne 0$ then $C^*$ has vanishing cohomology, as the complex is the
mapping cone of $g$ in degrees less than $\frac{n}{2}$ and the adjoint in
degrees greater than $\frac{n}{2}$. It follows that if $\epsilon \ne 0$ 
then the
Laplacian $d^C \left( d^C \right)^* + \left( d^C \right)^* d^C$ of 
$C^*$ has a bounded inverse.

Let $\nabla^W : W^* \rightarrow \Omega^1({\frak B}) \otimes_{\frak B} W^*$ 
be a self-dual connection on $W^*$. There is a direct sum connection
\begin{equation} \label{eqn75}
\nabla^C = \nabla^{{\cal D},0,1} \oplus \nabla^W
\end{equation}
on $C^*$.

Suppose that $n$ is even. Put $Q^C = d^C + \left( d^C \right)^*$.
We define a superconnection $D_s(\epsilon)$ on $C^*$ by
\begin{equation} \label{eqn76}
D_s(\epsilon) = s Q^C + \nabla^C.
\end{equation}

Let $\epsilon(s)$ be a smooth function of $s \in \R^+$ which is
identically zero for $s \in (0,1]$ and identically one for $s \ge 2$. Put
\begin{equation} \label{eqn77}
\widetilde{\eta}(s) = \TR_s \left( \frac{d D_s(\epsilon(s))}{ds} \:
e^{- D_s^2(\epsilon(s))} \right).
\end{equation}
\begin{proposition} \label{prop14}
For $s \in (0, 1]$, 
\begin{equation} \label{eqn78}
\widetilde{\eta}(s) = \TR_s \left( Q \: e^{-  D_s^2}
\right),
\end{equation}
as in (\ref{eqn69}).
\end{proposition}
\begin{pf}
As $\epsilon(s) = 0$, the factors $\Omega^*(F; {\cal D})$ and 
$\widehat{W}^*$ in
$C^*$ completely decouple and it is enough to show that the
analog of $\widetilde{\eta}(s)$ for $\widehat{W}^*$,
\begin{equation} \label{eqn79}
\TR_s \left( Q^{\widehat{W}} \: e^{-  (D_s^{\widehat{W}})^2}
\right),
\end{equation}
vanishes. This follows from
Hodge duality as in \cite[p. 227]{Lott (1992a)}. Namely, define
$T \in \End(\widehat{W}^*)$ to be multiplication by 
$\sign \left( i - \frac{n}{2} 
\right)$ on $\widehat{W}^i$. It is odd with respect to the Hodge duality
on $\widehat{W}^*$. Then 
\begin{align} \label{eqn80}
\TR_s \left( Q^{\widehat{W}} \: e^{-  (D_s^{\widehat{W}})^2} \right) & =
\TR_s \left( T^{-1} \: T \: Q^{\widehat{W}} \: 
e^{-  (D_s^{\widehat{W}})^2} 
\right)  = \:
- \: \TR_s \left(T \: Q^{\widehat{W}} \: e^{-  (D_s^{\widehat{W}})^2} 
T^{-1} \right)
\\
& = \: - \: \TR_s \left(T \: T^{-1} \: Q^{\widehat{W}} \: 
e^{-  (D_s^{\widehat{W}})^2} \right)  = \:
- \: \TR_s \left( Q^{\widehat{W}} \: e^{-  (D_s^{\widehat{W}})^2} 
\right)  = 0. 
\notag
\end{align}
\end{pf}

\begin{definition}
Define $\widetilde{\eta} \in \overline{\Omega}^{odd}({\frak B})
/\Image(d)$ by
\begin{equation} \label{eqn81}
\widetilde{\eta} = \int_0^\infty \widetilde{\eta}(s) \: ds. 
\end{equation}
\end{definition}
By \cite[Proposition 26]{Lott (1992a)}, 
the integrand of (\ref{eqn81}) is integrable for small-$s$.
Using the techniques of \cite[Section 6.1]{Lott (1996)} and
\cite[Section 4]{Lott (1998)}, one can show that
it is also integrable for large-$s$. This uses the 
invertibility of the Laplacian of $C^*$ for 
$s \ge 2$, i.e. $\epsilon = 1$.
Note that $\widetilde{\eta}$ is defined modulo $\Image(d)$.
It is not hard to show that $\widetilde{\eta}$ is independent of the
choice of $\epsilon(s)$.

\begin{proposition} \label{prop15}
$\widetilde{\eta}$ is independent of the choice of $W$.
\end{proposition}
\begin{pf}
Let $W^\prime$ be another regular Hermitian complex which is homotopy
equivalent to $\Omega^*(F; {\cal D})$, with $W^{\prime,\frac{n}{2}} = 0$. 
Let $h : W^\prime \rightarrow W$ be a homotopy equivalence.

For $-1 \le i \le n+1$, put 
\begin{equation} \label{eqn82}
D^i = 
\begin{cases}
\Omega^i(F; {\cal D}) \oplus W^i \oplus W^{i+1} \oplus W^{\prime,i+1} & 
\text{if $-1 \le i < \frac{n}{2}$,} \\
\Omega^{\frac{n}{2}}(F; {\cal D}) & \text{ if $n$ is even and $i =
\frac{n}{2}$,} \\
\Omega^i(F; {\cal D}) \oplus W^i \oplus W^{i-1} \oplus W^{\prime,i-1}
& \text{if $\frac{n}{2} < i \le n + 1$.} \\
\end{cases}
\end{equation}
Given 
$\begin{pmatrix}
\epsilon_1 & \epsilon_2 \\
\epsilon_3 & \epsilon_4
\end{pmatrix} \in M_2(\R)$, define a differential $d^D$ on $D^*$ by
\begin{equation} \label{eqn83}
d^D_i = 
\begin{cases}
\begin{pmatrix}
d & 0 & \epsilon_1 g & \epsilon_2 g \circ h \\
0 & d & \epsilon_3 & \epsilon_4 h \\
0 & 0 & -d & 0 \\
0 & 0 & 0 & -d
\end{pmatrix}
& \text{if $i < \frac{n}{2}$,} \\
d & \text{ if $n$ is even and $i = \frac{n}{2}$,} \\
\begin{pmatrix}
d & 0 & 0 & 0\\
0 & d & 0 & 0\\
\epsilon_1 g^* & \epsilon_3 & -d & 0 \\
\epsilon_2 h^* \circ g^*  & \epsilon_4 h^* & 0 & -d
\end{pmatrix}
& \text{if $i > \frac{n}{2}$.} \\
\end{cases}
\end{equation}
As
\begin{equation} \label{eqn84}
\begin{pmatrix}
\epsilon_1 g & \epsilon_2 g \circ h \\
\epsilon_3 & \epsilon_4 h
\end{pmatrix} =
\begin{pmatrix}
g & 0 \\
0 & 1
\end{pmatrix}
\begin{pmatrix}
\epsilon_1 & \epsilon_2  \\
\epsilon_3 & \epsilon_4 
\end{pmatrix}
\begin{pmatrix}
1 & 0 \\
0 &  h
\end{pmatrix},
\end{equation}
the complex $D^*$ has vanishing cohomology if 
$\begin{pmatrix}
\epsilon_1 & \epsilon_2  \\
\epsilon_3 & \epsilon_4 
\end{pmatrix}$ is invertible.
Given 
$A = \begin{pmatrix}
a_1 & a_2  \\
a_3 & a_4 
\end{pmatrix} \in \GL(2, \R)$, put
\begin{equation} \label{eqn85}
\begin{pmatrix}
\epsilon_1 & \epsilon_2  \\
\epsilon_3 & \epsilon_4 
\end{pmatrix} = \epsilon(s) \: A,
\end{equation}
where $\epsilon(s)$ is as before.
Define the noncommutative eta-form of $D^*$ as in (\ref{eqn81}).

There is a smooth path in $\GL(2, \R)$ from 
$\begin{pmatrix}
1 & 0  \\
0 & 1 
\end{pmatrix}$ to 
$\begin{pmatrix}
0 & 1  \\
-1 & 0 
\end{pmatrix}$.
It follows from \cite[(50)]{Lott (1992a)} that the corresponding eta-forms
of $D^*$ differ by something in $\Image(d)$. Consider the
eta-form coming from 
$A = \begin{pmatrix}
1 & 0  \\
0 & 1 
\end{pmatrix}$. In this case, $D^*$ splits into the sum of two complexes,
one involving $\Omega^*$ and $W^*$, the other involving $W^*$ and 
$W^{\prime,*}$. By the Hodge duality argument of Proposition \ref{prop14}, the
eta-form of the second subcomplex vanishes.  Hence when 
$A = \begin{pmatrix}
1 & 0  \\
0 & 1 
\end{pmatrix}$, we
recover the eta-form of the complex $C^*$ constructed from $\Omega^*$ and
$W^*$. Similarly, when $A = \begin{pmatrix}
0 & 1  \\
-1 & 0 
\end{pmatrix}$, we recover the eta-form constructed from $\Omega^*$ and
$W^{\prime,*}$. The proposition follows. 
\end{pf}

If $n$ is odd then one can define the higher eta-form using an extra
Clifford variable as in \cite[Definitions 2,11]{Lott (1992a)}.

\subsection{``Moral'' Fundamental Group of $S^1 \backslash M$}
\label{``Moral'' Fundamental Group of Quotient}

Let $M$ be a closed oriented smooth manifold with an effective $S^1$-action.
Let $\Gamma^\prime$ be a finitely generated discrete group and let 
$\rho : \pi_1(M) \rightarrow \Gamma^\prime$ be a surjective homomorphism.
There is an induced connected normal $\Gamma^\prime$-covering 
$M^\prime$ of $M$,
on which $\gamma^\prime 
\in \Gamma^\prime$ acts on the left by $L_{\gamma^\prime} \in 
\Diff(M^\prime)$.
Let $\pi : M^\prime \rightarrow M$ be the projection map.

Define a Lie group $G^\prime$ as in (\ref{eqn48}), with $G = S^1$.
As the generator of the $S^1$-action on $M$ can be lifted to a vector field
on $M^\prime$, there is a short exact sequence
\begin{equation} \label{eqn86}
1 \longrightarrow \Gamma^\prime \longrightarrow G^\prime \longrightarrow S^1 
\longrightarrow 1.
\end{equation}
The homotopy exact sequence of this fibration gives
\begin{equation} \label{eqn87}
1 \longrightarrow \pi_1 \left( G^\prime \right) \longrightarrow \Z 
\longrightarrow \Gamma^\prime \longrightarrow \pi_0
\left( G^\prime \right) \longrightarrow 1.
\end{equation}
Put $\widehat{\Gamma} = \pi_0 \left( G^\prime \right)$. We will think of
$\widehat{\Gamma}$ as the ``moral'' fundamental group of $S^1 \backslash M$,
although it may not be the same as $\pi_1(S^1 \backslash M)$; 
$\widehat{\Gamma}$ also appears in the work of Browder-Hsiang
\cite{Browder-Hsiang (1982)}. Fixing a
basepoint $m_0 \in M$, let $o$ be the homotopy class of the orbit of $m_0$ in
$\pi_1(M, m_0)$.
From (\ref{eqn87}), 
$\widehat{\Gamma} = \Gamma^\prime /\rho( \langle o \rangle)$, 
where
$\langle o \rangle$ is the central subgroup of $\pi_1(M, m_0)$ generated by
$o$. If $M^{S^1} \neq \emptyset$ then it is natural to take $m_0 \in M^{S^1}$,
showing that $\langle o \rangle = \{e\}$.

Let $G^\prime_0$ be the connected component of the identity of $G^\prime$.
It is a copy of either $S^1$ or $\R$. Put 
$\widehat{M} = G^\prime_0 \backslash M^\prime$. 
Then $\widehat{\Gamma}$ acts properly and cocompactly on $\widehat{M}$, with
$\widehat{\Gamma} \backslash \widehat{M} = S^1 \backslash M$. Let $p :
\widehat{M} \rightarrow S^1 \backslash M$ be the
quotient map.
Putting
$\widehat{M^{S^1}} = p^{-1} \left( M^{S^1} \right)$, we
can describe $\widehat{M^{S^1}}$ as the cover of $M^{S^1}$ induced from
the composite map 
\begin{equation} \label{eqn88}
\pi_1 \left( M^{S^1} \right) \longrightarrow \pi_1(M) 
\stackrel{\rho}{\longrightarrow} \Gamma^\prime \longrightarrow 
\widehat{\Gamma}.
\end{equation}
The complement $\widehat{M} - \widehat{M^{S^1}}$ has a natural orbifold
structure.

We construct certain differential forms on the strata of $S^1 \backslash M$.
Let $h \in C^\infty_0 \left( \widehat{M^{S^1}} \right)$ satisfy
\begin{equation} \label{eqn89}
\sum_{\widehat{\gamma} \in \widehat{\Gamma}} 
L_{\widehat{\gamma}}^* h = 1;
\end{equation}
it is easy to construct such functions. Let $N$ be a small neighborhood of
$M^{S^1}$ in $S^1 \backslash M$ which is diffeomorphic to the mapping
cylinder of a fiber bundle, whose fibers are weighted complex projective
spaces and whose base is $M^{S^1}$. Let $\widehat{N}$ be the preimage of
$N$ in $\widehat{M}$, with projection 
$\widehat{q} : \widehat{N} \rightarrow
\widehat{M^{S^1}}$. Consider $\widehat{q}^* h$ on $\widehat{N}$. It can be
extended to
a compactly-supported function $H$ on $\widehat{M}$ which is smooth, in
the orbifold sense, on $\widehat{M} - \widehat{M^{S^1}}$ and which
satisfies
\begin{equation} \label{eqn90}
\sum_{\widehat{\gamma} \in \widehat{\Gamma}} L_{\widehat{\gamma}}^* H = 1.
\end{equation}

Consider the group cochains
\begin{align} \label{eqn91} 
C^k ( \widehat{\Gamma} ) = 
\{ \tau : \: & \widehat{\Gamma}^{k+1} \rightarrow \R : 
\tau \: {\mbox{\rm is skew
and for all }} \left( \widehat{\gamma}_0, \ldots, \widehat{\gamma}_k \right) 
\in 
\widehat{\Gamma}^{k+1} \: 
{\mbox{\rm and }} z \in \widehat{\Gamma}, \\
& \tau \left( \widehat{\gamma}_0 z, \widehat{\gamma}_1 z, \ldots, 
\widehat{\gamma}_k z \right) = 
\tau \left( \widehat{\gamma}_0, \widehat{\gamma}_1,
\ldots, \widehat{\gamma}_k \right) \}. \notag
\end{align}
Suppose that $\tau$ is a cocycle, i.e.
\begin{equation} \label{eqn92}
\sum_{j=0}^{k+1}
(-1)^j \: \tau \left( \widehat{\gamma}_0, \ldots, 
\widehat{\widehat{\gamma}_j}, \ldots, 
\widehat{\gamma}_{k+1} \right) = 0.
\end{equation}

\begin{definition}
Define an orbifold form $\widehat{\omega} \in \Omega^k \left(
\widehat{M} - \widehat{M^{S^1}} \right)$ by
\begin{equation} \label{eqn93}
\widehat{\omega} = \sum_{\widehat{\gamma}_1, \ldots,\widehat{\gamma}_k}
\tau \left(\widehat{\gamma}_1, \ldots,\widehat{\gamma}_k, e \right)
\: L_{\widehat{\gamma}_1}^* dH \wedge \ldots \wedge
L_{\widehat{\gamma}_k}^* dH.    
\end{equation}  
Define
$\widehat{\mu} \in \Omega^k \left( \widehat{M^{S^1}} \right)$ by
\begin{equation} \label{eqn94}
\widehat{\mu} = \sum_{\widehat{\gamma}_1, \ldots,
\widehat{\gamma}_k}
\tau \left(\widehat{\gamma}_1, \ldots,\widehat{\gamma}_k, e \right)
\: L_{\widehat{\gamma}_1}^* dh \wedge \ldots \wedge
L_{\widehat{\gamma}_k}^* dh.    
\end{equation}  
\end{definition}
\begin{proposition} \label{prop16}
There are closed forms $\omega \in \Omega^k \left( (S^1 \backslash M) - M^{S^1}
\right)$ and $\mu \in \Omega^k \left( M^{S^1} \right)$ such
that $\widehat{\omega} = p^* \omega$ and $\widehat{\mu} = p^* \mu$.
\end{proposition}
\begin{pf}
The proof is as in \cite[Lemma 4]{Lott (1992)}. We omit the details.
\end{pf}

Let $N$ be a small neighborhood of $M^{S^1}$ in $S^1 \backslash M$ as above,
with projection $q : N \rightarrow M^{S^1}$. By construction,
$\omega \big|_N = q^* \mu$. By \cite[\S 7]{Sullivan (1977)},
the pair $\left( \omega, \mu \right)$ represents a class in 
$\HH^k(S^1 \backslash M; \R)$. We obtain a map $\phi : 
\HH^* \left( \widehat{\Gamma}; \R \right) \longrightarrow 
\HH^* \left( S^1 \backslash M; \R \right)$ given by
$\phi([\tau]) = [(\omega, \mu)]$.

\begin{proposition} \label{prop17} 
(Browder-Hsiang \cite[Theorem 1.1]{Browder-Hsiang (1982)})
There is a commutative diagram
\begin{equation} \label{eqn95}
\begin{array}{ccc}
\HH^* \left( \widehat{\Gamma}; \R \right) & \stackrel{\phi}{\longrightarrow} & 
\HH^* \left( S^1 \backslash M; \R \right) \\
\alpha \downarrow & & \beta \downarrow \\
\HH^* \left( \Gamma^\prime; \R \right) & \stackrel{\gamma}{\longrightarrow} & 
\HH^* \left(  M; \R \right)
\end{array}
\end{equation}
where the bottom
row of (\ref{eqn95}) 
comes from the map $M \rightarrow B \Gamma^\prime$ induced by
$\rho$, the left column of
(\ref{eqn95}) comes from the homomorphism $\Gamma^\prime \rightarrow 
\widehat{\Gamma}$ and
the right column of (\ref{eqn95}) is pullback.
\end{proposition}
\begin{pf}
Given $[\tau] \in \HH^k \left( \widehat{\Gamma}; \R \right)$, represent it
by a cocycle $\tau \in Z^k \left( \widehat{\Gamma}; \R \right)$. Let
$\alpha(\tau) \in Z^k \left( \Gamma^\prime; \R \right)$ be its pullback
to $\Gamma^\prime$.
Let $r : M^\prime \rightarrow \widehat{M}$ be the quotient map. Then
$(\beta \circ \phi) [\tau]$ is characterized by the $\Gamma^\prime$-invariant
closed form
\begin{equation} \label{eqn96}
r^* \widehat{\omega} = 
\sum_{\widehat{\gamma}_1, \ldots,\widehat{\gamma}_k \in \widehat{\Gamma}}
\tau \left(\widehat{\gamma}_1, \ldots,\widehat{\gamma}_k, e \right)
\: L_{\widehat{\gamma}_1}^* d r^* H \wedge \ldots  
L_{\widehat{\gamma}_k}^* d r^* H    
\end{equation}  
on $M^\prime$.
Let $K \in C^\infty_0 \left( M^\prime \right)$ satisfy
\begin{equation} \label{eqn97}
\sum_{g \in \rho( \langle o \rangle)} L_g^* K = r^* H.
\end{equation}
Let $\widetilde{\Gamma} \subset \Gamma^\prime$ be a set of 
representatives for the cosets $\rho(\langle o \rangle) \backslash 
\Gamma^\prime$.
Then
\begin{align} \label{eqn98}
r^* \widehat{\omega} & = 
\sum_{\widetilde{\gamma}_1, \ldots,\widetilde{\gamma}_k \in 
\widetilde{\Gamma}}
\alpha(\tau) \left(\widetilde{\gamma}_1, \ldots,\widetilde{\gamma}_k, e \right)
\: L_{\widetilde{\gamma}_1}^* dr^*H \wedge \ldots 
L_{\widetilde{\gamma}_k}^* dr^*H \\
& = 
\sum_{{\gamma^\prime}_1, \ldots, {\gamma^\prime}_k \in 
{\Gamma}^\prime}
\alpha(\tau) \left({\gamma^\prime}_1, \ldots,{\gamma^\prime}_k, e \right)
\: L_{{\gamma^\prime}_1}^* dK \wedge \ldots 
L_{{\gamma^\prime}_k}^* dK. \notag
\end{align}  
By \cite[Proposition 14]{Lott (1992)}, 
the last term in (\ref{eqn98}) is the lift of a
closed $k$-form on $M$ which represents $\gamma ( \alpha ([\tau]))$.
\end{pf}
\noindent
{\bf Remark :} Theorem 1.1 of \cite{Browder-Hsiang (1982)} is phrased in terms
of rational cohomology and is valid for any compact connected Lie group, not
just $S^1$. We expect that our proof could be extended to these cases. \\

If $\widehat{\Gamma}$ satisfies Assumption \ref{ass1}, 
let ${\cal D}$ be the canonical ${\frak B}$-vector bundle on $M^{S^1}$. 
As in (\ref{eqn59}), the function $h$ gives a partially-flat
connection on ${\cal D}$. Let $[\tau] \in \HH^k(\Gamma; \C)$ be 
represented by a cocycle $\tau \in Z^k(\Gamma; \C)$ as in
Assumption \ref{ass1}. Then the Chern form 
$\ch \left( \nabla^{\cal D} \right) \in \Omega^* \left( M^{S^1};
\overline{\Omega}^*({\frak B}) \right)$ satisfies
\begin{equation} \label{eqn99}
\mu \: = \: c(k) \: 
\langle \ch \left( \nabla^{\cal D} \right), Z_\tau \rangle
\end{equation}
for some nonzero $c(k) \in \R$ \cite[Proposition 3]{Lott (1996)}.

\subsection{Fixed-point-free Actions II}
\label{Fixed-point-free Actions II}

Suppose that the $S^1$-action has no fixed-points.  There is a 
$C^*_r \widehat{\Gamma}$-Hilbert module of orbifold differential forms on 
$\widehat{M}$ and a (tangential) signature operator, which has an index 
$\sigma_{S^1}(M) \in K_* \left(C^*_r \widehat{\Gamma} \right)$.
Suppose that $\widehat{\Gamma}$ satisfies Assumption \ref{ass1}. Then for
any $[\tau] \in \HH^* \left( \widehat{\Gamma}; \R \right)$, we can
consider the pairing $\langle \sigma_{S^1}(M), Z_\tau \rangle \in \R$.
\begin{proposition} \label{prop18}
Construct $\phi ([\tau]) \in \HH^*(S^1 \backslash M; \R)$ as in the
previous subsection. Given a suborbifold ${\cal O}$ of $S^1 \backslash M$
as in Subsection \ref{Fixed-point-free Actions}, 
let $\phi ([\tau]) \big|_{\cal O}$ denote the
pullback of $\phi [\tau]$ to ${\cal O}$. Then
\begin{equation} \label{eqn100}
\langle \sigma_{S^1}(M), Z_\tau \rangle = 
\sum_{\cal O} \frac{1}{m_{\cal O}}
\int_{\cal O} L({\cal O}) \cup \phi ([\tau]) \big|_{\cal O}.
\end{equation}
\end{proposition}
\begin{pf}
The method of proof is the same as in \cite{Lott (1992)}, which dealt
with the case when $\widehat{\Gamma}$ acts freely on a smooth 
$\widehat{M}$. The only difference is that the local analysis must now be done
on orbifolds, as in \cite{Kawasaki (1978)}. We omit the details.
\end{pf}
\noindent
{\bf Remark :} It seems likely that Proposition \ref{prop18} follows from
a general localization result and is true whenever $\widehat{\Gamma}$ satisfies
the Strong Novikov Conjecture; compare 
\cite[Theorem 2.6]{Rosenberg-Weinberger (1993)}.

\subsection{Semifree Actions II}
\label{Semifree Actions II}

Suppose that $S^1$ acts effectively and semifreely on $M$. 
If $F$ is a connected component of $M^{S^1}$, put 
\begin{equation} \label{eqn101}
\Gamma_F = \Image \left( \pi_1(F) \longrightarrow \pi_1(M) 
\longrightarrow \Gamma^\prime \longrightarrow \widehat{\Gamma} \right).
\end{equation}
Suppose that
$\Gamma_F$ satisfies Assumption \ref{ass1}, with smooth subalgebra
${\frak B}_F$ of $C^*_r \Gamma_F$, and that $F$ satisfies
Assumption \ref{ass2} with respect to $C^*_r \Gamma_F$. Construct
$\widetilde{\eta} \in \overline{\Omega}^*({\frak B}_F)/\Image(d)$ 
for the manifold
$F$ as in Subsection \ref{Higher Eta-Invariant}.  
\begin{definition}
Given $[\tau] \in \HH^k \left( \widehat{\Gamma}; \R \right)$,
represent it by a cocycle $\tau \in Z^k \left(\widehat{\Gamma}; \R \right)$. 
Construct $\omega_\tau \in \Omega^k \left( (S^1 \backslash M) - M^{S^1} 
\right)$ as in Proposition \ref{prop16}. Given a connected component $F$ of
$M^{S^1}$, let $\tau_F \in Z^k(\Gamma_F; \R)$ 
be the restriction of $\tau$. Suppose that the cyclic cocycle $Z_{\tau_F}$
extends to a cyclic cocycle on ${\frak B}_F$. Put
\begin{equation} \label{eqn102}
\langle \sigma_{S^1}(M), [\tau] \rangle =
\int_{S^1 \backslash M} L(T(S^1 \backslash M)) \wedge \omega_\tau
\: + \: c(k) \: \sum_F  \langle \widetilde{\eta}, Z_{\tau_F} \rangle \in \R. 
\end{equation}
We assume that $k \equiv \dim(M) -1 \mod{4}$ so that the 
integral in (\ref{eqn102}) can be nonzero.
\end{definition}

\begin{theorem} \label{prop19}
$\langle \sigma_{S^1}(M), [\tau] \rangle$ is independent of the choices of
$S^1$-invariant Riemannian metric on $M$ and function 
$H \in C^\infty_0 \left( \widehat{M} \right)$ on $\widehat{M}$.
\end{theorem}
\begin{pf}
The method of proof is that of Proposition \ref{prop11}.
Define $\mu_\tau \in \Omega^* \left( M^{S^1} \right)$ as in 
Proposition \ref{prop16}. Let $J_{S^1}(M)$ 
denote the first term in the right-hand-side of (\ref{eqn102}). We first show
the metric independence. Let 
$\{g(\epsilon)\}_{\epsilon \in [0,1]}$ be a smooth $1$-parameter family of
$S^1$-invariant metrics on $M$. 
For simplicity, assume that $M^{S^1}$ has one connected component $F$.
As in (\ref{eqn39}), let us write
\begin{align} \label{eqn103}
L \left( R^{TF}(\epsilon) + d\epsilon \wedge 
\partial_\epsilon \omega^{TF}  \right) & = 
a_1 + d\epsilon \wedge a_2, \\
L \left( \widehat{\Omega}_V(\epsilon) + d\epsilon \wedge 
\partial_\epsilon \widehat{\omega}_V  \right) 
& = b_1 + d\epsilon \wedge b_2. \notag
\end{align}
Then as in (\ref{eqn41}),
\begin{equation} \label{eqn104}
J_{S^1}(M) \big|_{\epsilon = 1} - J_{S^1}(M) \big|_{\epsilon = 0} =
- \int_0^1 d\epsilon \wedge
\int_F  a_2 \wedge \mu_\tau \wedge \int_Z b_1.
\end{equation}
The Atiyah-Singer families index theorem gives 
an equality in $\HH^{even}(F; \R)$:
\begin{equation} \label{eqn105}
\ch (\Ind (d + d^*)) = \int_Z b_1, 
\end{equation}
where $d + d^*$ denotes the family of vertical signature operators on 
the fiber bundle $S^1 \backslash SNF \rightarrow F$. \\ \\
{\bf Case I.} $\dim(M) - \dim(F) \equiv 2 \mod{4}$.\\

As $Z = \C P^{2N}$ for some
$N$, $\Ind (d+d^*)$ is a trivial complex line bundle on $F$. Then
\begin{equation} \label{eqn106}
J_{S^1}(M) \big|_{\epsilon = 1} - J_{S^1}(M) \big|_{\epsilon = 0} =
- \int_0^1 d\epsilon \wedge
\int_F  a_2 \wedge \mu_\tau.
\end{equation}
On the other hand, from \cite[Proposition 27]{Lott (1992a)},
\begin{equation} \label{eqn107}
c(k) \: \langle \widetilde{\eta}, Z_{\tau_F} \rangle \big|_{\epsilon = 1} \: - 
c(k) \: \langle \widetilde{\eta}, Z_{\tau_F} \rangle \big|_{\epsilon = 0} =
\int_0^1 d\epsilon \wedge
\int_F a_2 \wedge \mu_\tau.
\end{equation} 
The proposition follows in this case. \\ \\
{\bf Case II.} $\dim(M) - \dim(F) \equiv 0 \mod{4}$.\\

As $Z = \C P^{2N+1}$ for some $N$,
$\Ind (d+d^*) = 0$. Then
\begin{equation} \label{eqn108}
J_{S^1}(M) \big|_{\epsilon = 1} - J_{S^1}(M) \big|_{\epsilon = 0} = 0.
\end{equation} 
Equation (\ref{eqn107}) is again valid.
As $a_2$ is concentrated in degree congruent to $-1 \mod{4}$ and
$k \equiv \dim(F) - 1 \mod{4}$, we have $\int_F a_2 \wedge \mu_\tau = 0$.
The proposition follows in this case.

Now fix the metric and suppose that
$\{H(\epsilon)\}_{\epsilon \in [0,1]}$ is a smooth $1$-parameter family
of functions $H$ constructed as in (\ref{eqn90}). 
Construct the corresponding form $\mu_\tau \in \Omega^* \left(
[0,1] \times F \right)$.
Write
\begin{equation} \label{eqn109}
\mu_\tau = a_1 + d\epsilon \wedge a_2,
\end{equation}
where $a_1, a_2 \in \Omega^*(F)$ depend on $\epsilon$.
Then
\begin{equation} \label{eqn110}
J_{S^1}(M) \big|_{\epsilon = 1} - J_{S^1}(M) \big|_{\epsilon = 0} =
- \int_0^1 d\epsilon \wedge
\int_F  L \left( R^{TF} \right) \wedge a_2.
\end{equation}
From \cite[Proposition 27]{Lott (1992a)},
\begin{equation} \label{eqn111}
c(k) \: \langle \widetilde{\eta}, Z_{\tau_F} \rangle \big|_{\epsilon = 1} \: - 
c(k) \: \langle \widetilde{\eta}, Z_{\tau_F} \rangle \big|_{\epsilon = 0} =
\int_0^1 d\epsilon \wedge
\int_F  L \left( R^{TF} \right) \wedge a_2.
\end{equation}
The proposition follows. 
\end{pf}

\subsection{General $S^1$-Actions II}
\label{General $S^1$-Actions II}

Let $S^1$ act effectively on $M$. For each connected component $F$ of
$M^{S^1}$, define $\Gamma_F$ as in (\ref{eqn101}).
Suppose that $\Gamma_F$ satisfies 
Assumption \ref{ass1}, with smooth subalgebra ${\frak B}_F \subset C^*_r 
\Gamma_F$, and that $F$ satisfies Assumption \ref{ass2} with respect to
$C^*_r \Gamma_F$. 
\begin{definition}
Given $[\tau] \in \HH^k \left( \widehat{\Gamma}; \R \right)$,
represent it by a cocycle $\tau \in Z^k \left(\widehat{\Gamma}; \R \right)$. 
Construct $\omega_\tau \in \Omega^k \left( (S^1 \backslash M) - M^{S^1} 
\right)$ as in Proposition \ref{prop16}. Given a connected component $F$ of
$M^{S^1}$, let $\tau_F \in Z^k(\Gamma_F; \R)$ 
be the restriction of $\tau$. Suppose that the cyclic cocycle $Z_{\tau_F}$
extends to a cyclic cocycle on ${\frak B}_F$. Put
\begin{equation} \label{eqn112}
\langle \sigma_{S^1}(M), [\tau] \rangle =
\sum_{\cal O} \frac{1}{m_{\cal O}} \int_{\cal O} L({\cal O}) \wedge
\omega_\tau \big|_{\cal O} 
\: + \: \sum_F \langle \widetilde{\eta}, Z_{\tau_F} \rangle \in \R. 
\end{equation}
\end{definition}

As in Theorem \ref{prop19}, $\langle \sigma_{S^1}(M), [\tau] \rangle$ is
independent of the choices of $S^1$-invariant metric and $H$.

\begin{conjecture}
$\langle \sigma_{S^1}(M), [\tau] \rangle$ is an  $S^1$-homotopy invariant
of $M$.
\end{conjecture}

One may want to assume that $\widehat{\Gamma}$ satisfies Assumption \ref{ass1}.
In this case,
if the $S^1$-action has no fixed-points then the conjecture follows from
Proposition \ref{prop18}, along with the homotopy invariance of the index
$\sigma_{S^1}(M) \in K_* \left( C^*_r \widehat{\Gamma} \right)$. If
the $S^1$-action is semifree and the codimension of $M^{S^1}$
in $M$ is at most two
then an outline of a proof of the conjecture is given in Appendix A.

\section{Remarks} \label{Remarks}
\noindent
1. One may wonder whether Assumption \ref{ass2} is really necessary.  To see
that some assumption is necessary to define equivariant higher signatures, 
consider the special case when the quotient space is a manifold-with-boundary.
So consider compact oriented manifolds-with-boundary equipped with a map to a
classifying space $B\pi$.
As Shmuel Weinberger pointed out to me, if one had a reasonable
higher signature for such manifolds then one would expect to have
Novikov additivity for the higher signatures of closed oriented 
manifolds. That is,
if $M$ is a closed oriented manifold with a map to $B\pi$ and $N$ is a 
hypersurface
in $M$ which cuts it into two pieces $M_1$ and $M_2$ then the higher signatures
of $M$ would be the sum of those of $M_1$ and $M_2$, for the same reasons that
the Atiyah-Patodi-Singer theorem implies the Novikov additivity of the usual
signature.  In particular, the higher signatures of closed oriented manifolds 
would give invariants of the cut-and-paste group $SK_*(B\pi)$
\cite{Karras-Kreck-Neumann-Ossa (1973)}. However, it is known
for some groups $\pi$ that the only cut-and-paste invariants of $B\pi$ are the
Euler characteristic and the usual signature. For example, it easy to
show that this is the case when $\pi = \Z$ and it then follows from
\cite[Lemma 8]{Neumann (1975)} that it is also the case when $\pi = \Z^k$.
Thus in general
one needs some assumption in order to define the higher signatures.

As a side remark, in some cases it is possible to define higher signatures of
manifolds-with-boundary without any extra assumptions.
For example, let $M$ be a compact oriented manifold-with-boundary such that
$4 | \dim(M)$. Let $\nu : M \rightarrow B\pi$ be a continuous map. 
Suppose that we are given a homomorphism
$\rho : \pi \rightarrow SO(p,q)$ for some $p, q > 0$.
Let $BSO(p,q)_\delta$ be the classifying space for $SO(p,q)$ with the
discrete topology. There is a canonical flat real vector bundle $V$ on
$BSO(p,q)_\delta$ of rank $p + q$. The pullback $(B\rho \circ \nu)^* V$ 
is a flat 
real vector bundle on $M$ with a flat symmetric form. Hence one can consider
the twisted signature $\sigma(M, (B\rho \circ \nu)^* V) \in \Z$. 
This is an oriented-homotopy
invariant of $M$ by construction. On the other hand, 
if $M$ is closed then it is also
a higher signature of $M$ involving the pullback of a Borel class from
$BSO(\infty, \infty)_\delta$ \cite{Lusztig (1971)}. It follows from the
usual Novikov additivity argument that this higher signature is a
cut-and-paste invariant.  For example, if $M$ is closed,
$4 | \dim(M)$ and
$B\pi$ is a closed oriented
hyperbolic manifold of dimension $\dim(M)$ then one finds
that the degree of the map $\nu$ gives a nontrivial invariant of 
$SK_{dim(B\pi)} B \pi$.
If $\dim(M) \equiv
2 \mod{4}$ then one can do a similar construction in which $SO(p,q)$ is
replaced by $Sp(2n)$. In general, it seems to be an interesting question
as to which higher signatures of closed manifolds are cut-and-paste
invariants.\\ \\
2. Although we have defined the signature of an $S^1$-quotient, we have not
defined a signature operator of which the signature is the index. If
$M^{S^1}$ has codimension in $M$ divisible by four then there is a signature
operator on $S^1 \backslash M$ by the work of Cheeger
\cite{Cheeger (1983)}. If the $S^1$-action is semifree and $M^{S^1}$ has
codimension two in $M$ then $S^1 \backslash M$ is a manifold-with-boundary
and one has the Atiyah-Patodi-Singer signature operator on $S^1 \backslash M$.
For a general semifree $S^1$-action, the quotient space will contain families
of cones over complex projective spaces.  We note that there is a topological
obstruction to having a self-adjoint signature operator on a singular
space with a single cone over $\C P^N$, $N$ even \cite{Lesch (1993)}.
However, in our case such cones occur in odd-dimensional families and
this fact may allow one to construct the signature operator.\\ \\
3. Suppose that a compact Lie group $G$ acts effectively on an oriented
closed manifold $M$. Let $M^{sing}$ be the set of points in $M$ whose
isotropy subgroup has positive dimension. Then we can define
$\Omega^{*,basic} \left( M, M^{sing} \right)$ and $\HH^{*,basic} 
\left( M, M^{sing} \right)$ as in Definition \ref{def1}. There is again an
intersection form on $\HH^{*,basic} \left( M, M^{sing} \right)$
which comes from integrating on the orbifold
$(G \backslash M) - (G \backslash M^{sing})$,
and its signature $\sigma_G(M)$ is a $G$-homotopy invariant of $M$.
One can ask for an explicit formula for $\sigma_G(M)$, as was done in this
paper when $G = S^1$. If the $G$-action is semifree then the analog of
Theorem \ref{prop7} holds and the proof is virtually the same as that of
Theorem \ref{prop7}. However, if the action is not semifree then the
situation is more involved.  Suppose, for simplicity, that all isotropy
groups are connected.  In principle, one can follow the proof of Theorem
\ref{prop7} 
by applying the Atiyah-Patodi-Singer formula to a sequence of compact 
manifolds-with-boundary that exhaust
$(G \backslash M) - (G \backslash M^{sing})$. However, the limiting
formula must be more complicated than in Theorem \ref{prop7}. For example
take $G = SU(2)$. If $m \in M^{sing} - M^{SU(2)}$ then a neighborhood of
$\overline{m} \in SU(2) \backslash M$ is like an $S^1$-quotient of the type
studied in Section \ref{Signatures of $S^1$-Quotients}. 
In analogy to Theorem \ref{prop7}, we expect that
there will be a contribution to $\sigma_G(M)$ of the form
$\eta \left( (G \backslash M^{sing}) - (G \backslash M^{SU(2)})\right)$. 
However, $(G \backslash M^{sing}) - (G \backslash M^{SU(2)})$ is a space
with conical singularities like those in Section 
\ref{Signatures of $S^1$-Quotients} and it is not 
immediately clear how to define its eta-invariant; this is related to the 
preceding remark.
\appendix
\section{Homotopy Invariance of Higher Signatures of Manifolds-With-Boundary}
\label{Homotopy Invariance of Higher Signatures of Manifolds-With-Boundary}

Suppose that we have a compact oriented 
manifold-with-boundary $A$, a finitely generated discrete group
$\widehat{\Gamma}$ and a surjective homomorphism 
$\pi_1(A) \rightarrow \widehat{\Gamma}$. For simplicity, suppose that
$\partial A$ just has one connected component. Put 
\begin{equation} \label{eqn113}
\Gamma_F = \Image \left( \pi_1(\partial A) \longrightarrow \pi_1(A) 
\longrightarrow \widehat{\Gamma} \right).
\end{equation}
Put $n = \dim(\partial A)$.
Suppose that $\widehat{\Gamma}$ satisfies Assumption \ref{ass1}, with
smooth subalgebra ${\frak B}$ of $C^*_r \widehat{\Gamma}$. Put
\begin{equation}
{\frak B}_F = \{ T \in {\frak B} : T \left( l^2(\Gamma_F) \right) \subset
l^2(\Gamma_F)\}.
\end{equation}
Then $\C \Gamma_F \subset {\frak B}_F \subset C^*_r \Gamma_F$, with
${\frak B}_F$ closed under the holomorphic
functional calculus in $C^*_r \Gamma_F$. Let $i : \overline{\Omega}^* 
({\frak B}_F)/\Image(d) \rightarrow 
\overline{\Omega}^* ({\frak B})/\Image(d)$ be the obvious map.
Suppose that $\partial A$ satisfies Assumption \ref{ass2} with respect to
$C^*_r \Gamma_F$. Construct $\widetilde{\eta} \in \overline{\Omega}^* 
({\frak B}_F)/\Image(d)$ for 
$\partial A$ as in Subsection \ref{Higher Eta-Invariant}.
Let ${\cal D}$ be the canonical flat ${\frak B}$-vector bundle on $A$.
We have the higher signature
\begin{equation} \label{eqn114}
\sigma(A) = 
\int_{A} L(TA) \wedge \ch \left(
\nabla^{\cal D} \right)
\: + \: i(\widetilde{\eta}) \: \in \: \overline{\HH}^*({\frak B}).
\end{equation}

We want to realize $\sigma(A)$ as the Chern character of an index.
We first describe the ``unperturbed'' setting.  Without loss of generality,
suppose that $A$ is metrically a product near $\partial A$. 
Put $B = A \cup_{\partial A} \left( [0, \infty) \times \partial A \right)$.
We extend ${\cal D}$ over $B$ as a product over the cylindrical end. 
Consider the ${\frak B}$-module 
$\Omega^* \left( B; {\cal D} \right)$ of smooth compactly-supported 
${\cal D}$-valued forms on $B$. This
is one component of the unperturbed situation. 

We would like to interpret
$\sigma(A)$ as the index of the
signature operator on the $C^*_r \widehat{\Gamma}$-completion of
$\Omega^* \left( B; {\cal D} \right)$. However, there is the problem that
this signature operator need not be Fredholm in the 
$C^*_r \widehat{\Gamma}$-sense, because the signature operator on
$\Omega^* \left( \partial A; {\cal D} \right)$ may not be invertible.
This is why we proceed as follows.

The other component of the unperturbed situation
is an algebraic analog of a half-infinite cylinder which is coned off.
More precisely, let $W^*$ be a cochain complex of finitely generated
projective ${\frak B}$-modules which
is homotopy equivalent to $\Omega^*(\partial A; {\cal D})$ as in 
Subsection \ref{Higher Eta-Invariant}. Let $\widehat{W}^*$ be as in 
(\ref{eqn70}).
Let $\phi \in C^\infty([0, \infty))$ be a nondecreasing
function such that
\begin{equation} \label{eqn115}
\phi(r) = 
\begin{cases}
r \text{ if $r \le \frac{1}{2}$,} \\
1 \text{ if $r \ge 2$.} \\ 
\end{cases}
\end{equation}
Define a ${\frak B}$-inner product on the ${\frak B}$-cochain complex
$\Omega^*((0, \infty)) {\otimes} \widehat{W}^*$ such that if
$w^i \in C^\infty_0((0, \infty)) \otimes \widehat{W}^i$ and 
$w^j \in C^\infty_0((0, \infty)) \otimes \widehat{W}^j$ then
\begin{align} \label{eqn116}
\langle w^i, w^j \rangle & = \int_0^\infty \phi(r)^{dim(\partial A) - i - j} 
\: 
\langle w^i(r), w^j(r) \rangle_{\widehat{W}} \: dr, \\
\langle w^i, dr \: \wedge \: w^j \rangle & = 0, \notag \\
\langle dr \: \wedge \: w^i, dr \: \wedge \: w^j \rangle & = 
\int_0^\infty \phi(r)^{dim(\partial A) - i - j} \: 
\langle w^i(r), w^j(r) \rangle_{\widehat{W}} \: dr. \notag
\end{align}
This is the second component of the unperturbed situation.
Formally, one would expect from Hodge duality that the index of a signature
operator on $\Omega^*((0, \infty)) {\otimes} \widehat{W}^*$ 
should vanish. Hence the index of a signature operator on 
$\Omega^* \left( B; {\cal D} \right) \oplus
\left( \Omega^*((0, \infty)) {\otimes} \widehat{W}^* \right)$
is formally the same as that of
$\Omega^* \left( B; {\cal D} \right)$.

We now perturb $\Omega^* \left( B; {\cal D} \right) \oplus
\left( \Omega^*((0, \infty)) {\otimes} \widehat{W}^* \right)$ 
to obtain a Fredholm operator.  Let $d$ denote the total differential on
$\Omega^* \left( B; {\cal D} \right) \oplus
\left( \Omega^*((0, \infty)) {\otimes} \widehat{W}^* \right)$,
where we switch the sign on the $\widehat{W}^*$-differential as in 
(\ref{eqn72}).
Let $\epsilon \in C^\infty([0, \infty))$ be a nondecreasing
function such that
\begin{equation} \label{eqn117}
\epsilon(r) = 
\begin{cases}
0 \text{ if $r \le 1$,} \\
1 \text{ if $r \ge 2$.} \\ 
\end{cases}
\end{equation}
Given $\alpha > 1$,
define an operator $D$ on 
$\Omega^* \left( B; {\cal D} \right) \oplus
\left( \Omega^*((0, \infty)) {\otimes} \widehat{W}^* \right)$ by
saying that on the degree-$i$ subspace,
\begin{equation} \label{eqn118}
D = d + 
\begin{cases}
\begin{pmatrix}
0 & \epsilon(r/\alpha) \:  g \\
0 & 0
\end{pmatrix}
& \text{if $i < \frac{n}{2}$,} \\
& \\
\begin{pmatrix}
0 & 0 \\
\epsilon(r/\alpha) \: f & 0
\end{pmatrix}
& \text{if $i > \frac{n}{2}$.} \\
\end{cases}
\end{equation}
Note that $D^2 \ne 0$ because $\epsilon$ is a nonconstant function of $r$.
If $n + 1$ is even, put $T = D + D^*$. If $n + 1$ is odd, 
put $T = \pm (*D - D*)$.
Then we expect that it will be possible to show the following :\\
1. The operator $T$ extends to a Fredholm operator in the 
$C^*_r \widehat{\Gamma}$-sense. Its index $\Ind(T)$ 
is independent of $\alpha$.\\
2. In analogy to \cite{Leichtnam-Piazza (1997)}, 
$\ch(\Ind(T)) = \sigma(A)$.\\
3. In analogy to \cite{Hilsum-Skandalis (1992)}, $\Ind(T)$ is a smooth
homotopy invariant of the pair $(A, \partial A)$. (That is, the 
homotopy equivalence is not required to be a diffeomorphism on $\partial A$.)
\\ 

To relate this to $S^1$-actions, let $M$ have a semifree $S^1$-action such 
that $M^{S^1}$ is nonempty and has 
codimension two.  Then $S^1 \backslash M$ is a manifold-with-boundary
$A$, with $\partial A = M^{S^1}$. 
By \cite[Proposition 1.2]{Rosenberg-Weinberger (1988)},
$\pi_1(M) = \pi_1(A)$. If $\rho : \pi_1(M) \rightarrow \Gamma^\prime$ is
a surjective homomorphism as in Section 
\ref{``Moral'' Fundamental Group of Quotient} 
then $\widehat{\Gamma} = \Gamma^\prime$.

Extending point 3. above, 
we mean that $\Ind(T)$ should be an $S^1$-homotopy invariant of
$M$. Given an $S^1$-homotopy equivalence $h : M \rightarrow N$, put
$A = S^1 \backslash M$ and $B = S^1 \backslash N$. We obtain
a homotopy equivalence 
$\overline{h} : A \rightarrow B$ on the quotient spaces. 
It may not be a proper map, in that  
$\partial A = M^{S^1}$ may be properly contained
in the preimage of $\partial B = N^{S^1}$.
Nevertheless, we can extend $\overline{h}$ to a smooth map 
$h^\prime : A \cup_{\partial A} ([0, \infty) \times \partial A) \rightarrow
B \cup_{\partial B} ([0, \infty) \times \partial B)$ which is a product map
on $[0, \infty) \times \partial A$. It should be possible
to use $h^\prime$, as in \cite{Hilsum-Skandalis (1992)}, to compare the
signature operators of $A$ and $B$. The analog of the almost-flat connection
of \cite{Hilsum-Skandalis (1992)} is the fact that although $D^2 \ne 0$,
by taking $\alpha$ large we can make the norm of $D^2$ as small as we want.
Regarding point 3. above, 
it may be more convenient to work with a conical end
than a cylindrical end.  This
would correspond to multiplying the metric on $[0, \infty) \times \partial A$
by a conformal factor which is asymptotically $e^{-2cr}$ for large $r$, and
similarly changing the inner product on 
$\left( \Omega^*((0, \infty)) {\otimes} \widehat{W}^* \right)$.
Here $c$ is some positive constant.
\\ \\
{\bf Remark :} In the topological setting, with similar assumptions one
has a symmetric signature $s(A) \in L^* \left(\Z \widehat{\Gamma} \right)$.
To describe this, assume for simplicity that $\Gamma_F = \widehat{\Gamma}$.
Following \cite{Weinberger (199?)}, assume that $\partial A$ is antisimple,
meaning that the chain complex 
$C_* \left( \partial A; \Z \widehat{\Gamma} \right)$ is 
homotopy equivalent to a chain
complex $P_*$ of finitely generated projective $\Z \widehat{\Gamma}$-modules,
with $P_{\frac{n}{2}} = 0$ if $n$ is even and 
$P_{\frac{n \pm 1}{2}} = 0$ if $n$ is odd.  
Let $P_<$ denote the truncation of $P_*$ at
$\left[ \frac{n}{2} \right]$. Then the map 
$C_* \left( \partial A; \Z \widehat{\Gamma} \right) \rightarrow P_<$ defines
an algebraic Poincar\'e pair in the sense of
\cite[p. 134]{Ranicki (1980)}. The (closed)
algebraic Poincar\'e complex $C_* \left(A;  \Z \widehat{\Gamma} \right) 
\cup_{C_* \left( \partial A; \Z \widehat{\Gamma} \right)} P_<$ has a
symmetric signature $s(A) \in L^* \left(\Z \widehat{\Gamma} \right)$ which
will be a homotopy invariant of the pair $(A, \partial A)$. 
If $\widehat{\Gamma}$ satisfies Assumption \ref{ass1} then we can construct
$\ch (s(A)) \in \overline{\HH}_*({\frak B})$; compare with (\ref{eqn114}).

\end{document}